\documentclass[12pt,a4paper]{article}

\usepackage{amssymb}

\usepackage{amsmath}

\newtheorem{theorem}{Theorem}[section]
\newtheorem{definition}[theorem]{Definition}
\newtheorem{lemma}[theorem]{Lemma}
\newtheorem{corollary}[theorem]{Corollary}

\newtheorem{example}[theorem]{Example}
\newtheorem{conjecture}[theorem]{Conjecture}

\begin{document}
\title{Gr\"{o}bner-Shirshov bases for Lie algebras over a commutative algebra\footnote{Supported by the
NNSF of China (Nos.10771077, 10911120389) and the NSF of Guangdong
Province (No. 06025062).}}
\author{
L.A. Bokut\footnote {Supported by RFBR 01-09-00157, LSS--344.2008.1
and SB RAS Integration grant No. 2009.97 (Russia).} \\
{\small \ School of Mathematical Sciences, South China Normal
University}\\
{\small Guangzhou 510631, P. R. China}\\
{\small Sobolev Institute of Mathematics, Russian Academy of
Sciences}\\
{\small Siberian Branch, Novosibirsk 630090, Russia}\\
{\small Email: bokut@math.nsc.ru}\\
\\
Yuqun Chen\footnote {Corresponding author.} \ and Yongshan Chen
\\
{\small \ School of Mathematical Sciences, South China Normal
University}\\
{\small Guangzhou 510631, P. R. China}\\
{\small Email: yqchen@scnu.edu.cn}\\
{\small jackalshan@126.com}}

\date{}

\maketitle \noindent\textbf{Abstract:} In this paper we establish a
Gr\"{o}bner-Shirshov bases theory for Lie algebras over  commutative
rings. As applications we give some new examples of special Lie
algebras (those embeddable in associative algebras over the same
ring) and non-special Lie algebras (following a suggestion of P.M.
Cohn (1963) \cite{Conh}). In particular, Cohn's Lie algebras over
the characteristic $p$ are non-special when $p=2,\ 3,\ 5$. We
present an algorithm that one can check for any $p$, whether Cohn's
Lie algebras is non-special. Also we prove that any finitely or
countably generated Lie algebra is embeddable in a two-generated Lie
algebra.

\noindent \textbf{Key words: } Lie algebra over a commutative ring,
Lyndon-Shirshov word, Gr\"{o}bner-Shirshov basis.

\noindent \textbf{AMS 2000 Subject Classification}: 17B01, 16S15,
13P10

\section{Introduction}

Gr\"{o}bner bases and Gr\"{o}bner-Shirshov bases  were invented
independently by A.I. Shirshov \cite{Sh62a,Shir3} for ideals of free
(commutative, anti-commutative) non-associative algebras, free Lie
algebras  \cite{Sh62b,Shir3} and implicitly free associative
algebras  \cite{Sh62b,Shir3}  (see also \cite{Be78,Bo76}), by H.
Hironaka \cite{Hi64} for ideals of the power series algebras (both
formal and convergent), and by B. Buchberger \cite{Bu70} for ideals
of the polynomial algebras.

 The Shirshov's Composition-Diamond lemma and Buchberger's theorem is
the corner stone of the theories. This proposition says that in
appropriate free algebra $A_{\bf k}(X)$ over a field ${\bf k}$ with
a free generating set $X$ and a fixed monomial ordering, the
following conditions on a subset $S$ of $A_{\bf k}(X)$ are
equivalent:
\begin{enumerate}
\item[(i)]\  Any
composition ($s$-polynomial) of  polynomials from $S$ is trivial;
\item[(ii)]\  If $f\in Id(S)$, then the maximal monomial $\bar{f}$ contains
some maximal monomial $\bar s$, where $s\in S$ (for Lie algebra
case, $\bar f$ means the maximal associative word of Lie polynomial
$f$);
\item[(iii)]\  The set
$Irr(S)$ of all (non-associative in general) words in $X $, which do
not contain any maximal word $\bar s, s\in S$, is a linear $k$-basis
of the algebra $A(X|S)=A(X)/Id(S)$ with generators $X$ and defining
relations $S$ (for Lie algebra case, $Irr(S)$ is the set of
Lyndon--Shirshov Lie words whose associative supports do not contain
maximal associative words of polynomials from $S$).
\end{enumerate}

$S$ is called a Gr¡§obner-Shirshov basis of the ideal $Id(S)$ of
$A_k(X)$ generated by $S$ if one of the conditions (i)-(iii) holds.

Gr\"{o}bner bases and Gr\"{o}bner-Shirshov bases theories have been
proved to be very useful in different branches of mathematics,
including commutative algebra and combinatorial algebra, see, for
example, the books
 \cite{AL, BKu94, BuCL, BuW, CLO, Ei}, the papers \cite{Be78, Bo72,Bo76},
 and the surveys \cite{BC, BFKK00, BK03, BK05}.

Up to now, different versions of Composition-Diamond lemma are known
for the following classes of algebras apart those mentioned above:
(color) Lie super-algebras (\cite{Mik89, Mik92,Mik96}), Lie
$p$-algebras \cite{Mik92}, associative conformal algebras
\cite{BFK}, modules \cite{KL,Chi} (see also \cite{CCZ},
right-symmetric algebras \cite{BCLi08}, dialgebras \cite{BCL08},
associative algebras with multiple operators \cite{BCQ08},
Rota-Baxter algebras \cite{BCD08}, and so on.

It is well-known Shirshov's result \cite{Shir1, Shir3} that every
finitely or countably generated Lie algebra over a field ${\bf k}$
can be embedded into a two-generated Lie algebra over ${\bf k}$.
 Actually, from the technical
point of view, it was a beginning of the Gr\"{o}bner-Shirshov bases
theory for Lie algebras (and associative algebras as well). Another
proof of the result using explicitly Gr\"{o}bner-Shirshov bases
theory is refereed
 to  L.A. Bokut, Yuqun Chen  and
Qiuhui Mo \cite{BCM}.

A.A. Mikhalev and  A.A. Zolotykh \cite{MZ} prove the
Composition-Diamond lemma for a tensor product of a free algebra and
a polynomial algebra, i.e., they establish Gr\"{o}bner-Shirshov
bases theory for
 associative algebras over
 a commutative algebra. L.A. Bokut, Yuqun Chen and Yongshan Chen \cite{BCC08} prove the
Composition-Diamond lemma for a tensor product of two free algebras.
Yuqun Chen, Jing Li and Mingjun Zeng \cite{CJZ} prove the
Composition-Diamond lemma for a tensor product of a non-associative
algebra and a polynomial algebra.

In this paper, we establish the Composition-Diamond lemma for Lie
algebras over a polynomial algebra, i.e., for ``double free" Lie
algebras. It provides a Gr\"{o}bner-Shirshov bases theory for Lie
algebras over a commutative algebra.

Let ${\bf k}$ be a field, $K$ a commutative associative ${\bf
k}$-algebra with identity, and $\mathcal{L}$ a Lie $K$-algebra. Let
$Lie_K(X)$ be the free Lie $K$-algebra generated by a set $X$. Then,
of course, $\mathcal{L}$ can be presented as $K$-algebra by
generators $X$ and some defining relations $S$,
$$
                              \mathcal{L}=Lie_K(X|S)=Lie_K(X)/Id(S).
$$
In order to define a Gr\"{o}bner-Shirshov basis for $\mathcal{L}$, we first present $K$ in a form

$$
                           K={\bf k}[Y|R]={\bf k}[Y]/Id(R),
$$
where ${\bf k}[Y]$ is a polynomial algebra  over the field ${\bf k}$, $R\subset{\bf k}[Y]$. Then the
Lie $K$-algebra $\mathcal{L}$ has the following presentation as a ${\bf k}[Y]$-algebra
$$
          \mathcal{L}= Lie_{{\bf k}[Y]}(X| S, Rx,\ x\in X)
$$
(cf.  E.S. Chibrikov \cite{Chi}, see also \cite{CCZ}).

Now by definition, a Gr\"{o}bner-Shirshov basis for $\mathcal{L}=
Lie_K(X|S)$ is Gr\"{o}bner-Shirshov basis (in the sense of the
present paper) of the ideal $Id(S, Rx,\ x\in X)$
 in the ``double free" Lie algebra $Lie_{{\bf k}[Y]}(X)$.

As an application of our Composition-Diamond lemma (Theorem \ref{cdL2}), a  Gr\"{o}bner-Shirshov basis
of $\mathcal{L}$ gives rise to a linear basis of $\mathcal{L}$ as a ${\bf k}$-algebra.

We give applications of Gr\"{o}bner-Shirshov bases theory for Lie
algebras over a commutative algebra $K$ (over a field ${\bf k}$) to
the Poincare-Birkhoff-Witt theorem. Recent survey on PBW theorem see
in P.-P. Grivel \cite{Grivel}.  A Lie algebra over a commutative
ring is called special if it is embeddable into an (universal
enveloping)  associative algebra.  Otherwise it is called
non-special. There are known  classical examples by A.I. Shirshov
\cite{Shir53} and P. Cartier  \cite{Cartier} of Lie algebras over
commutative algebras over $GF(2)$ that are not embeddable into
associative algebras. Shirshov and Cartier  used ad hoc methods to
prove that some elements of corresponding  Lie algebras are not zero
though they are zero in the universal enveloping algebras, i.e.,
they proved non-speciality of the examples. Here we find
Gr\"{o}bner-Shirshov bases of these Lie algebras and then use our
Composition-Diamond lemma to get the result, i.e.,  we give a new
conceptual proof.

P.M. Cohn \cite{Conh} gave the following examples of Lie algebras
$$
\mathcal{L}_p=Lie_{K}(x_1,x_2,x_3|y_3x_3=y_2x_2+y_1x_1)
$$
over truncated polynomial algebras
$$
K={\bf k}[y_1,y_2,y_3 | y_i^p=0, 1\leq i\leq 3],
$$
where ${\bf k}$ is a filed of characteristic $p>0$. He conjectured
that $\mathcal{L}_p$ is non-special Lie algebra for any $p$.
$\mathcal{L}_p$ is called the Cohn's Lie algebra. Using our
Composition-Diamond lemma we have proved that $\mathcal{L}_2$,
$\mathcal{L}_3$ and $\mathcal{L}_5$ are non-special Lie algebras. We
present an algorithm that one can check for any $p$, whether Cohn's
Lie algebras is non-special.

We give new class of special Lie algebras in terms of defining
relations (Theorem \ref{t4.5}). For example, any one relator Lie
algebra $Lie_K(X| f)$ with a ${\bf k}[Y]$-monic relation $f$ over a
commutative algebra $K$ is special (Corollary \ref{co4.6}). It gives
an extension of the list of known special Lie algebras (ones with
valid PBW Theorems) (see P.-P. Grivel \cite{Grivel}). Let us give
this list:
\begin{enumerate}
\item[1.]\
  $\mathcal{L}$ is a free $K$-module
(G. Birkhoff \cite{Birkhoff}, E. Witt \cite{Witt}),
\item[2.]\   $K$ is a principal ideal domain (M. Lazard \cite{Lazard52, Lazard54}),
\item[3.]\    $K$ is a Dedekind domain  (P. Cartier \cite{Cartier}),
\item[4.]\  $K$ is over a field ${\bf k}$ of characteristic 0 (P.M. Cohn \cite{Conh}),
\item[5.]\   $\mathcal{L}$ is $K$-module without torsion  (P.M. Cohn \cite{Conh}),
\item[6.]\  $2$ is invertible in $K$ and for any $x,y,z\in \mathcal{L}$,
$[x[yz]]=0$ (Y. Nouaze and P. Revoy \cite{NouazeRevoy}).
\end{enumerate}

P. Higgins \cite{Higgins} unified the cases 1-3 and gave homological
invariants of special Lie algebras inspired by results of R. Baer,
see also P. Revoy \cite{Revoy}.

As a last application we prove that every finitely or countably
generated Lie algebra over an arbitrary commutative algebra $K$  can
be embedded into a two-generated Lie algebra over $K$.

We thank Yu Li and Jiapeng Huang for some comments.

\section{Preliminaries}

We start with some concepts and results from the literature
concerning with the Gr\"{o}bner-Shirshov bases theory of a free Lie
algebra $Lie_{{\bf k}}(X)$ generated by $X$ over a field ${\bf k}$.

Let $X=\{x_i |i\in I\}$ be a well-ordered set with $x_i>x_j$ if
$i>j$ for any $i,j\in I$. Let $X^*$ be the free monoid generated by
$X$. For $u=x_{i_1}x_{i_2} \cdots x_{i_m} \in X^*$, let the  length
of $u$ be $m$, denoted by $|u|=m.$

 We use two linear orderings on $X^*$:

(i) ({\sl lex ordering}) $1>t$ if $t\neq1$ and, by induction, if
$u=x_{i}u^{\prime}$ and $v=x_{j}v^{\prime}$ then $u>v$ if and only
if $x_{i}>x_{j}$ or $x_{i}=x_{j}$ and $u^{\prime}>v^{\prime}$;

(ii) ({\sl deg-lex ordering}) $u \succ v$ if $|u|>|v|$, or $|u|=|v|$
and $u>v$.\\

We regard $Lie_{{\bf k}}(X)$ as the Lie subalgebra of the free
associative algebra ${\bf k} \langle X \rangle$, which is generated
by $X$ under the Lie bracket $[u,v]=uv-vu$. Given $f \in {\bf k}
\langle X \rangle$, denote by $\bar{f}$ the leading word of $f$ with
respect to the deg-lex ordering; $f$ is  {\sl monic} if the
coefficient of $\bar{f}$ is 1.

\begin{definition}(\cite{ Lyndon, Shir1})
$w\in X^*\setminus\{1\}$ is an associative Lyndon--Shirshov
 word (ALSW for short) if $$(\forall u,v\in X^*, u,v\neq1) \ w=uv\Rightarrow w>vu.$$
\end{definition}
We denote the set of all ALSW's on $X$ by $ALSW(X)$.

We cite some useful properties of ALSW's (\cite{Lyndon,Shir1}, see
also, for example, \cite{bc07, BK03, BK05, BKu94, Reutenauer, U95}):

({\bf I}) if $w\in ALSW(X)$   then an arbitrary proper prefix of $w$
cannot be a suffix of $w$;

({\bf II}) if $w=uv\in ALSW(X)$, where $u,v \ne 1$ then $u>w>v$;

({\bf III}) if $u,v\in ALSW(X)$ and $u>v$ then $uv\in ALSW(X)$;

({\bf IV}) an arbitrary associative word $w$ can be uniquely
represented as $w=c_{1}c_{2} \ldots c_{n},$ where $c_{1}, \ldots
,c_{n}\in ALSW(X)$ and $c_{1} \leq c_{2} \leq \ldots \leq c_{n}$;

({\bf V}) if $u'=u_1u_2$ and $u''=u_2u_3$ are ALSW's then
$u=u_1u_2u_3$ is also an ALSW;

({\bf VI}) if an associative word $w$ is represented as in ({\bf
IV}) and $v$ is an associative Lyndon-Shirshov  subword of $w$, then
$v$ is a subword of one of the words $c_{1}$,
$c_{2}$,$\ldots$,$c_{n}$;

({\bf VII}) if an ALSW $w=uv$ and $v$ is its  longest proper ALSW,
then $u$ is an ALSW as well.

\begin{definition} (\cite{CFL, Shir1})\
A non-associative word $(u)$ in $X$ is a non-associative
Lyndon-Shirshov word (NLSW for short), denoted by $[u]$,  if

(i) $u$ is an ALSW;

(ii) if $[u]=[(u_{1})(u_{2})]$ then both $(u_{1})$ and $(u_{2})$ are
NLSW's (from ({\bf I}) it then follows that $u_{1}>u_{2}$);

(iii) if $[u]=[[[u_{11}][u_{12}]][u_{2}]]$ then $u_{12} \leq u_{2}$.
\end{definition}
We denote the set of all NLSW's on $X$ by $NLSW(X)$.

In fact,  NLSW's may be defined as Hall--Shirshov words relative to
lex ordering (for definition of Hall--Shirshov words see
\cite{Sh62c}, also \cite{VG}).

By \cite{Lyndon, Shir1, Shir3}, for an ALSW $w$, there is a unique
bracketing $[w]$ such that $[w]$ is NLSW:  $[w]=w$ if $|w|=1$ and
$[w]=[[u][v]]$ if $|w|>1$, where $v$ is the longest proper
associative Lyndon-Shirshov end of $w$ and by $({\bf VII}) $ $u$ is
an ALSW. Then by induction on $|w|$, we have $[w]$.

It is well-known that  the set $NLSW(X)$ forms a linear basis of
$Lie_{{\bf k}}(X)$, see  \cite{Lyndon, Shir1, Shir3}.

Considering any NLSW $[w]$ as a polynomial in ${\bf k} \langle X
\rangle$, we have $\overline{[w]}=w$ (see \cite{Shir1, Shir3}).
 This implies that  if $f \in Lie_{{\bf k}}(X)\subset {\bf k} \langle X \rangle$
then $\bar{f}$  is an ALSW.

\begin{lemma}(Shirshov \cite{Shir1, Shir3})\label{l1} Suppose that $w=aub$,
where $w,u\in ALSW(X)$. Then
$$
[w]=[a[uc]d],
$$
where  $b=cd$ and possibly $c=1$. Represent $c$ in the form
$$
c=c_{1}c_{2} \ldots c_{n},
$$
where $c_{1}, \ldots ,c_{n}\in ALSW(X)$ and $c_{1} \leq c_{2} \leq
\ldots \leq c_{n}$. Replacing $[uc]$ by $[\ldots[[u][c_{1}]] \ldots
[c_{n}]]$ we obtain the word $[w]_{u}=[a[\ldots[[[u][c_{1}]][c_2]]
\ldots [c_{n}]]d]$ which is called the special bracketing of $w$
 relative to  $u$. We have
$$
\overline{[w]}_{u}=w.
$$
\end{lemma}

\begin{lemma}(Chibrikov \cite{Chi2})\label{chi'lem1} Let $w=aub$ be as in Lemma \ref{l1}.
Then $[uc]=[u[c_1][c_2]\ldots[c_n]]$, that is
$$
[w]=[a[\ldots[u[c_{1}]] \ldots
[c_{n}]]d].
$$
\end{lemma}

\begin{lemma}(\cite{BKu94,Chi2})\label{chi'lem} Suppose that $w=aubvc$, where $w,u,v\in ALSW(X)$. Then
there is some bracketing
$$
[w]_{u,v}=[a[u]b[v]d]
$$
in the word $w$ such that
$$
\overline{[w]}_{u,v}=w.
$$
More precisely,
$$
[w]_{u,v}=\left\{
\begin{array} {ll}\ [a[up]_uq[vs]_vl] &\mbox{ if }
[w]=[a[up]q[vs]l],\\
\ [a[u[c_1]\cdots[c_t]_v\cdots[c_n]]_up] &\mbox{ if }   [w]=
[a[u[c_1]\cdots[c_t]\cdots[c_n]]p] \mbox{ with }  v  \mbox{ a
subword of } c_t.
\end{array}\right.
$$
\end{lemma}

\section{Composition-Diamond lemma for $Lie_{{\bf k}[Y]}(X)$}

Let $Y=\{y_j|j\in
J\}$ be a well-ordered set and $[Y]=\{y_{j_1}y_{j_2}\cdots
y_{j_l}|y_{j_1}\leq y_{j_2}\leq\cdots \leq y_{j_l}, l \geq0
\}$ the free commutative monoid generated by $Y$. Then $[Y]$ is a
${\bf k}$-linear basis of the polynomial algebra ${\bf k}[Y]$.

Let the set $X$ be a well-ordered set, and let the lex ordering $<$ and the deg-lex ordering
$\prec_X$
 on $X^*$ be defined as before.

Let $Lie_{{\bf k}[Y]}(X)$ be the ``double" free Lie algebra, i.e., the free Lie algebra over the polynomial algebra ${\bf k}[Y]$ with generating set $X$.

From now on we regard $Lie_{{\bf k}[Y]}(X)\cong {\bf k}[Y]\otimes Lie_{{\bf k}}(X)$ as
the Lie subalgebra of ${\bf k}[Y] \langle X \rangle\cong{\bf k}[Y]\otimes{\bf k} \langle X \rangle$ the free
associative algebra over polynomial algebra ${\bf k}[Y]$, which is
generated by $X$ under the Lie bracket $[u,v]=uv-vu$.

Let
$$
T_A=\{u=u^Yu^X|u^Y\in[Y], \ u^X  \in ALSW(X)\}
$$
and
$$
T_N=\{[u]=u^Y[u^X]|u^Y\in[Y], \ [u^X]\in NLSW (X)\}.
$$

By the previous section, we know that the elements of $T_A$ and $T_N$ are one-to-one corresponding to each other.

\noindent {\bf Remark:} For $u=u^Yu^X\in T_A$, we still use the
notation $[u]=u^Y[u^X]$ where $[u^X]$ is a NLSW on $X$.

Let ${\bf k}T_N$ be the linear space spanned by $T_N$ over ${\bf
k}$. For any $[u],[v]\in T_N$, define
$$
[u][v]=\sum\alpha_iu^Yv^Y[w_i^X]
$$
where $\alpha_i\in {\bf k},\ [w_i^X]$'s are NLSW's and
$[u^X][v^X]=\sum\alpha_i[w_i^X]$ in $Lie_{{\bf k}}(X)$.

Then ${\bf k}[Y]\otimes Lie_{{\bf k}}(X)\cong {\bf k}T_N$ as ${\bf k}$-algebra
and $T_N$ is a ${\bf k}$-basis of ${\bf k}[Y]\otimes Lie_{{\bf k}}(X)$.

We define the deg-lex ordering $\succ$ on
$$[Y]X^*=\{u^Yu^X|u^Y\in
[Y],u^X\in X^*\}
$$
by the following: for any $u,v\in [Y]X^*$,
\begin{eqnarray*}
&&u\succ v  \ \mbox{ if } \ (u^X\succ_Xv^X)  \  \mbox{ or }\
(u^X=v^X \  \mbox{ and }\ u^Y\succ_Yv^Y ),
\end{eqnarray*}
where $\succ_Y$ and $\succ_X$ are the deg-lex ordering on $[Y]$ and $X^*$ respectively.

\noindent {\bf Remark:} By abuse of the notation, from now on, in a
Lie expression like
 $[[u][v]]$ we will omit the external brackets, $[[u][v]]=[u][v]$.

Clearly, the ordering $\succ$ is ``monomial" in a sense of
$\overline{[u][w]}\succ\overline{[v][w]}$ whenever $w^X\neq u^X$ for
any $u,v,w\in T_A$.

Considering any $[u]\in T_N$ as a polynomial in ${\bf k}$-algebra
${\bf k}[Y]\langle X\rangle$, we have
$\overline{[u]}=u\in T_A$.

For any $f\in Lie_{{\bf k}[Y]}(X)\subset{\bf k}[Y]\otimes {\bf
k}\langle X\rangle$, one can present $f$ as a ${\bf k}$-linear
combination of $T_N$-words, i.e.,  $f=\sum\alpha_i[u_i]$, where
$[u_i]\in T_N$. With respect to the ordering $\succ$ on $[Y]X^*$,
the leading word $\bar{f}$ of $f$ in ${\bf k}[Y]\langle X\rangle$ is
an element of $T_A$. We call $f$  ${\bf k}$-monic if the coefficient
of $\bar{f}$ is 1. On the other hand, $f$ can be presented as ${\bf
k}[Y]$-linear combinations of $NLSW(X)$, i.e., $f=\sum
f_i(Y)[u_i^X]$, where $f_i(Y)\in {\bf k}[Y]$, $[u_i^X]\in NLSW(X)$
and $u_1^X \succ_Xu_2^X\succ_X\ldots$. Clearly $\bar{f}^X=u_1^X$ and
$\bar{f}^Y=\overline{f_1(Y)}$. We call $f$  ${\bf k}[Y]$-monic if
the $f_1(Y)=1$. It is easy to see that ${\bf k}[Y]$-monic implies
${\bf k}$-monic.

Equipping with the above concepts, we rewrite Lemma \ref{l1} as
follows.

\begin{lemma}\label{l2}(Shirshov \cite{Shir1, Shir3}) Suppose that $w=aub$ where $w,u\in T_A$ and $a,b\in X^*$. Then
$$
[w]=[a[uc]d],
$$
where $[uc]\in T_N$  and $b=cd$.

Represent $c$ in a form $c=c_{1}c_{2} \ldots c_{n},$ where $c_{1},
\ldots ,c_{n}\in ALSW(X)$ and $c_{1} \leq c_{2} \leq \ldots \leq
c_{n}$. Then
$$
[w]=[a[u[c_{1}][c_{2}] \ldots [c_{n}]]d].
$$
Moreover, the leading word of $[w]_u=[a[\cdots[[[u][c_{1}]][c_{2}] ]\ldots [c_{n}]]d]$
is exactly $w$, i.e.,
$$
\overline{[w]}_{u}=w.
$$
\end{lemma}

We still use the notion $[w]_u$ as the special bracketing of $w$
 relative to  $u$ in Section $2$.

Let $S\subset Lie_{{\bf k}[Y]}(X) $ and $Id(S)$ be the ${\bf
k}[Y]$-ideal of $Lie_{{\bf k}[Y]}(X)$ generated by $S$. Then
any element of $Id(S)$ is a ${\bf k}[Y]$-linear combination of polynomials of the
following form:
$$
(u)_{s}=[c_1][c_2]\cdots[c_{n}]s[d_1][d_2]\cdots[d_{m}], \ \
m,n\geq0
$$
with some placement of parentheses, where $s\in S$ and $c_i,d_j\in
ALSW(X)$. We call such $(u)_s$ an $s$-word (or $S$-word).

Now, we define two special kinds of $S$-words.
\begin{definition} Let $S\subset Lie_{{\bf k}[Y]}(X)$ be a ${\bf k}$-monic subset, $a,b\in X^*$ and
$s\in S$. If $a\bar{s}b\in T_A$,  then by Lemma \ref{l2} we have the
special bracketing $[a\bar{s}b]_{\bar{s}}$ of $a\bar{s}b$ relative
to  $\bar{s}$.  We define
$[asb]_{\bar{s}}=[a\bar{s}b]_{\bar{s}}|_{[\bar{s}]\mapsto{s}}$ to be
a normal $s$-word (or normal $S$-word).
\end{definition}

\begin{definition} Let $S\subset Lie_{{\bf k}[Y]}(X)$ be a ${\bf k}$-monic
subset and $s\in S$. We define the quasi-normal $s$-word, denoted by
$\lfloor u\rfloor_{{s}}$, where $u=asb, \ a,b\in X^*$ ($u$ is an
associative $S$-word), inductively.
\begin{enumerate}
\item[(i)] $s$ is quasi-normal of $s$-length $1$;
\item[(ii)] If $\lfloor u\rfloor_s$ is quasi-normal with $s$-length $k$ and
$[v]\in NLSW(X)$ such that $|v|=l$, then  $[v]\lfloor u\rfloor_s$
when $v>\overline{\lfloor u\rfloor}_s^X$ and $\lfloor
u\rfloor_s[v]$ when $v<\overline{\lfloor u\rfloor}_s^X$ are
quasi-normal of $s$-length $k+l$.
 \end{enumerate}
\end{definition}

From the definition of the quasi-normal $s$-word, we have the
following lemma.

\begin{lemma}\label{l7} For any quasi-normal $s$-word $\lfloor u\rfloor_s=(asb), \ a,b\in X^*$, we have
 $\overline{\lfloor u\rfloor}_s=a\bar{s}b\in T_A$.
\end{lemma}

\noindent {\bf Remark:} It is clear that for an $s$-word
$(u)_{s}=[c_1][c_2]\cdots[c_{n}]s[d_1][d_2]\cdots[d_{m}]$, $(u)_{s}$
is quasi-normal if and only if $\overline{(u)_s}=c_1c_2\cdots
c_{n}\overline{s}d_1d_2\cdots d_{m}$.

\ \

Now we give the definition of compositions.

\begin{definition}\label{d1}
Let $f,g$ be  two ${\bf k}$-monic polynomials of $Lie_{{\bf
k}[Y]}(X)$. Denote the least common multiple of $\bar{f}^Y$ and
$\bar{g}^Y$ in $[Y]$ by $L=lcm(\bar{f}^Y,\bar{g}^Y)$.

If $\bar{g}^X$ is a subword of $ \bar{f}^X$, i.e.,
$\bar{f}^X=a\bar{g}^Xb$ for some $a,b\in X^*$, then the polynomial
$$
C_1\langle f,g\rangle_w=\frac{L}{\bar{f}^Y}f-\frac{L}{\bar{g}^Y}[agb]_{\bar{g}}
$$
is called the\textbf{ inclusion composition} of $f$ and $g$ with
respect to $w$, where $w=L\bar{f}^X=La\bar{g}^Xb$.

If a proper prefix of $\bar{g}^X$ is a proper suffix of $\bar{f}^X$,
i.e., $\bar{f}^X=aa_0$, $\bar{g}^X=a_0b$, $a,b,a_0\neq1$, then the
polynomial
$$
C_2\langle f,g\rangle_w=\frac{L}{\bar{f}^Y}[fb]_{\bar{f}}-\frac{L}{\bar{g}^Y}[ag]_{\bar{g}}
$$
is called the \textbf{intersection composition} of $f$ and $g$ with
respect to $w$, where $w=L\bar{f}^Xb=La\bar{g}^X$.

 If the greatest common divisor of $\bar{f}^Y$ and $\bar{g}^Y$ in
$[Y]$ is not 1, then for any $a,b,c\in X^*$ such that
$w=La\bar{f}^Xb\bar{g}^Xc\in T_A$, the polynomial
$$
C_3\langle f,g\rangle_w=\frac{L}{\bar{f}^Y}[afb\bar{g}^Xc]_{\bar{f}}-
\frac{L}{\bar{g}^Y}[a\bar{f}^Xbgc]_{\bar{g}}
$$
is called the \textbf{external composition} of $f$ and $g$ with
respect to $w$.

 If  $\bar{f}^Y\neq1$,
then for any normal $f$-word $[afb]_{\bar{f}}, \ a,b\in X^*$, the polynomial
$$
C_4\langle f\rangle_w=[a\bar{f}^Xb][afb]_{\bar{f}}
$$
is called the \textbf{multiplication composition} of $f$ with
respect to $w$, where $w=a\bar{f}^Xba\bar{f}b$.
\end{definition}

Immediately, we have that $\overline{C_i\langle -\rangle_w}\prec w,
\ i\in \{1,2,3,4\}.$

 \ \

\noindent{\bf Remarks:}\begin{enumerate}
\item[1)]
 When $Y=\varnothing$, there are no external
and multiplication compositions. This is the case of Shirshov's
compositions over a field.

\item[2)] In the cases of $C_1$ and $C_2$, the corresponding $w\in T_A$ by the property of ALSW's,
but in the case of $C_4$, $w\not\in T_A$.

\item[3)] For any fixed $f,g$, there are finitely many compositions
$C_1\langle f,g\rangle_w$, $C_2\langle f,g\rangle_w$, but infinitely
many $C_3\langle f,g\rangle_w$, $C_4\langle f\rangle_w$.

\end{enumerate}

\begin{definition}
Given a ${\bf k}$-monic subset $S\subset Lie_{{\bf k}[Y]}(X)$ and
$w\in[Y]X^* $ (not necessary in $T_A$), an element $h\in Lie_{{\bf
k}[Y]}(X)$ is called trivial modulo $(S,w)$, denoted by $h\equiv0 \
mod(S,w)$, if $h$ can be presented as a ${\bf k}[Y]$-linear
combination of normal $S$-words with leading words less than $w$,
i.e., $h=\sum_{i}\alpha_{i}\beta_i[a_{i}s_{i}b_{i}]_{\bar{s}_{i}}$,
where $\alpha_{i} \in {\bf k}$, $\beta_i\in [Y]$, $a_{i}, b_{i} \in
X^*$, $s_{i} \in S$, and $\beta_ia_{i}\bar{s}_{i}b_{i}\prec w$.

In general, for $p,q\in Lie_{{\bf k}[Y]}(X)$, we write $p\equiv q \ mod(S,w)$ if $p-q\equiv0 \ mod(S,w)$.

$S$ is a Gr\"{o}bner-Shirshov basis in $Lie_{{\bf k}[Y]}(X)$ if all the possible compositions of elements in $S$ are
trivial modulo $S$ and corresponding $w$.
\end{definition}

If a subset $S$ of $Lie_{{\bf k}[Y]}(X)$ is not a
Gr\"{o}bner-Shirshov basis then one can add all nontrivial
compositions of polynomials of $S$ to $S$. Continuing this process
repeatedly, we finally obtain a Gr\"{o}bner-Shirshov basis $S^{C}$
that contains $S$. Such a process is called Shirshov's algorithm.
$S^{C}$ is called Gr\"{o}bner-Shirshov complement
 of $S$.

\begin{lemma}\label{l40} Let $f$ be a ${\bf k}$-monic polynomial in $Lie_{{\bf k}[Y]}(X)$.
If $\bar{f}^Y=1$ or $f=gf'$ where $g\in {\bf k}[Y]$ and
$f'\in Lie_{{\bf k}}(X)$, then for any normal $f$-word
$[afb]_{\bar{f}}$, $a,b\in X^*$,
$(u)_f=[a\bar{f}^Xb][afb]_{\bar{f}}$ has a presentation:
$$
(u)_f=[a\bar{f}^Xb][afb]_{\bar{f}}=\sum_{\overline{\lfloor
u_i\rfloor}_{f}\preceq\overline{(u)}_f}\alpha_i\beta_i\lfloor
u_i\rfloor_{f}
$$
where $\alpha_{i} \in {\bf k}, \ \beta_i\in [Y].$
\end{lemma}
{\bf Proof.} Case 1. $\bar{f}^Y=1$, i.e., $\bar{f}=\bar{f}^X$. By
Lemma \ref{l2} and since $\prec$ is monomial, we have
$[a\bar{f}b]=[afb]_{\bar{f}}-\sum_{\beta_iv_i\prec
a\bar{f}b}\alpha_i\beta_i[v_i]$, where $\alpha_i\in {\bf k}, \
\beta_i\in[Y], \ v_i\in ALSW(X)$. Then
$$
(u)_f=[a\bar{f}b][afb]_{\bar{f}}\\
=[afb]_{\bar{f}}[afb]_{\bar{f}}+\sum_{\beta_iv_i\prec a\bar{f}b }\alpha_i\beta_i[afb]_{\bar{f}}[v_i]\\
=\sum_{\beta_iv_i\prec a\bar{f}b}\alpha_i\beta_i[afb]_{\bar{f}}[v_i].
$$
The result follows since $v_i\prec a\bar{f}b$ and each $[afb]_{\bar{f}}[v_i]$ is quasi-normal.

Case 2. $f=gf'$, i.e., $\bar{f}^X=\bar{f'}$.  Then we have
$$
(u)_f=[a\bar{f'}b][afb]_{\bar{f}} =g([a\bar{f'}b][af'b]_{\bar{f'}}).
$$
The result follows from Case 1. \ \ \ $\Box$

\ \

The following lemma plays a key role in this paper.

\begin{lemma}\label{l4}Let $S$ be a ${\bf k}$-monic subset of $Lie_{{\bf k}[Y]}(X)$ in which
each multiplication composition is trivial. Then for any
quasi-normal $s$-word $\lfloor u\rfloor_s=(asb)$ and $w=a\bar{s}b=\overline{\lfloor u\rfloor_s}$,
where $a,b\in X^*$, we have
$$
\lfloor u\rfloor_s \equiv[asb]_{\bar{s}}  \ \  mod(S,w).
$$

\end{lemma}

{\bf Proof.} For $w=\bar s$  the lemma is clear.

For $w\neq\bar s$, since either $\lfloor u\rfloor_s=(asb)=[a_1](a_2sb)$ or
$\lfloor u\rfloor_s=(asb)=(asb_1)[b_2]$, there are
two cases to consider.

Let
$$
\delta_{(asb)}=\left\{
\begin{array} {ll} |a_1| &\mbox{if } (asb)=[a_1](a_2sb),\\
s\mbox{-length of } (asb_1) &\mbox{if } (asb)=(asb_1)[b_2].
\end{array} \right.
$$
The proof will be proceeding  by induction on $(w,\delta_{(asb)})$,
where   $(w',m')<(w,m)\Leftrightarrow
w\prec w'$ or $w=w',\ m'<m$ ($w,w'\in T_A,m,m'\in \mathbb{N}$).

Case 1. $\lfloor u\rfloor_s=(asb)=[a_1](a_2sb)$, where
$a_1>a_2\overline{s}^Xb,\ a=a_1a_2$ and $(a_2sb)$ is quasi-normal
$s$-word. In this case, $(w,\delta_{(asb)})=(w,|a_1|)$.

Since $w=a\bar{s}b=a_1a_2\bar{s}b\succ a_2\bar{s}b$, by
induction, we may assume that $
(a_2sb)=[a_2sb]_{\bar{s}}+\sum\alpha_i\beta_i[c_is_id_i]_{\bar{s_i}}
$, where $\beta_ic_i\bar{s_i}d_i\prec a_2\bar{s}b, \ a_1,a_2,c_i,d_i\in X^*, \ \ s_i\in S,\
\alpha_i\in {\bf k}$ and $\beta_i\in [Y]$. Thus,
$$
\lfloor
u\rfloor_s=(asb)=[a_1][a_2sb]_{\bar{s}}+\sum\alpha_i\beta_i[a_1][c_is_id_i]_{\bar{s_i}}.
$$

Consider the term $[a_1][c_is_id_i]_{\bar{s_i}}$.

If $a_1>c_i\bar{s_i}^Xd_i$, then $[a_1][c_is_id_i]_{\bar{s_i}}$ is
quasi-normal $s$-word with $a_1c_i\bar{s_i}d_i\prec w$. Note that
$\beta_ia_1c_i\bar{s_i}d_i\prec w$, then by induction,
$\beta_i[a_1][c_is_id_i]_{\bar{s_i}}\equiv0 \ mod(S,w)$.

If $a_1<c_i\bar{s_i}^Xd_i$, then
$[a_1][c_is_id_i]_{\bar{s}_i}=-[c_is_id_i]_{\bar{s}_i}[a_1]$ and
$[c_is_id_i]_{\bar{s_i}}[a_1]$ is quasi-normal $s$-word with
$\beta_ic_i\bar{s_i}d_ia_1\prec \beta_ia_2\bar{s}ba_1\prec\beta_i
a_1a_2\bar{s}b=w$.

If $a_1=c_i\bar{s_i}^Xd_i$, then there are two possibilities. For
${s_i}^Y=1$, by Lemma \ref{l40} and  by induction on $w$ we have
$\beta_i[a_1][c_is_id_i]_{\bar{s_i}}\equiv0 \ \ mod(S,w)$. For
${s_i}^Y\neq1$, $[a_1][c_is_id_i]_{\bar{s_i}}$ is the multiplication
composition, then by assumption, it is trivial $mod(S,w)$.

This shows that in any case, $\beta_i[a_1][c_is_id_i]_{\bar{s_i}}$
is a linear combination of normal $s$-words with leading words less
than $w$, i.e., $\beta_i[a_1][c_is_id_i]_{\bar{s_i}}\equiv0 \ \
mod(S,w)$ for all $i$.

Therefore, we may assume that $\lfloor
u\rfloor_s=(asb)=[a_1][a_2sb]_{\bar{s}}$ and $a_1>w^X>
a_2\bar{s}^Xb$.

If either $|a_1|=1$ or $[a_1]=[[a_{11}][a_{12}]]$ and $a_{12}\leq
a_2\bar{s}^Xb$, then $\lfloor u\rfloor_s=[a_1][a_2sb]_{\bar{s}}$ is
already a normal $s$-word, i.e., $\lfloor
u\rfloor_s=[a_1][a_2sb]_{\bar{s}}=[a_1a_2sb]_{\bar{s}}=[asb]_{\bar{s}}$.

If $[a_1]=[[a_{11}][a_{12}]]$ and $a_{12}>a_2\bar{s}^Xb$, then
$$
\lfloor u\rfloor_s=[a_1][a_2sb]_{\bar{s}}=[[a_{11}][a_{12}]][a_2sb]_{\bar{s}}=[a_{11}][[a_{12}][a_2sb]_{\bar{s}}]
+[[a_{11}][a_2sb]_{\bar{s}}][a_{12}].
$$

Let us consider the second summand
$[[a_{11}][a_2sb]_{\bar{s}}][a_{12}]$. Then by induction on $w$ and
by noting that $[a_{11}][a_2sb]_{\bar{s}}$ is quasi-normal, we may
assume that
$[a_{11}][a_2sb]_{\bar{s}}=\sum\alpha_i\beta_i[c_is_id_i]_{\bar{s_i}}
$, where $\beta_ic_i\bar{s_i}d_i\preceq a_{11}a_2\overline{s}b,\
s_i\in S,\ \alpha_i\in {\bf k},  \ \beta_i\in [Y], \ c_i,d_i\in
X^*$. Thus,
$$
[[a_{11}][a_2sb]_{\bar{s}}][a_{12}]=\sum\alpha_i\beta_i[c_is_id_i]_{\bar{s_i}}[a_{12}],
$$
where $ a_{11}>a_{12}>a_2\bar{s}^Xb,\
w=a_{11}a_{12}a_2\overline{s}b.$

If $a_{12}<c_i\bar{s_i}^Xd_i$, then
$[c_is_id_i]_{\bar{s_i}}[a_{12}]$ is quasi-normal with
$w'=\beta_ic_i\bar{s}_id_ia_{12}\preceq\beta_i
a_{11}a_2\bar{s}ba_{12}\prec w$. By induction,
$\beta_i[c_is_id_i]_{\bar{s_i}}[a_{12}]\equiv0 \ mod(S,w).$

If $a_{12}>c_i\bar{s_i}^Xd_i$, then
$[c_is_id_i]_{\bar{s_i}}[a_{12}]=-[a_{12}][c_is_id_i]_{\bar{s_i}}$
and $[a_{12}][c_is_id_i]_{\bar{s_i}}$ is quasi-normal with
$w'=\beta_ia_{12}c_i\bar{s}_id_i\preceq
\beta_ia_{12}a_{11}a_2\bar{s}b\prec w$. Again we can apply the
induction.

If $a_{12}=c_i\bar{s_i}^Xd_i$, then as discussed above, it is either
the case in Lemma \ref{l40} or the multiplication composition and
each is trivial $mod(S,w)$.

These show that $[[a_{11}][a_2sb]_{\bar{s}}][a_{12}]\equiv 0 \ \
mod(S,w)$.

Hence,
$$
\lfloor u\rfloor_s\equiv
 [a_{11}][[a_{12}][a_2sb]_{\bar{s}}] \ \ \ mod(S,w).
$$
where $a_{11}>a_{12}>a_2\bar{s}^Xb.$

Noting that $[a_{11}][[a_{12}][a_2sb]_{\bar{s}}]$ is quasi-normal and now
$(w,\delta_{[a_{11}][[a_{12}][a_2sb]_{\bar{s}}]})=(w,|a_{11}|)<(w,|a_1|)$, the result follows  by
induction.

\vspace{0.5cm}

Case 2. $\lfloor u\rfloor_s=(asb)= (asb_1)[b_2]$ where
$a\overline{s}^Xb_1>b_2,\ b=b_1b_2$ and $(asb_1)$ is quasi-normal
$s$-word. In this case, $(w,\delta_{(asb)})=(w,m)$ where $m$ is the
$s$-length of $(asb_1)$.

By induction on $w$, we may assume that
$$
\lfloor
u\rfloor_s=(asb)=[asb_1]_{\bar{s}}[b_2]+\sum\alpha_i\beta_i[c_is_id_i]_{\overline{s_i}}[b_2].
$$
where $\beta_ic_i\overline{s_i}d_i\prec a\overline{s}b_1,\ s_i\in S,\
\alpha_i\in {\bf k}, \ \beta_i\in[Y], \ c_i,d_i\in X^*$.

Consider the term $\beta_i[c_is_id_i]_{\overline{s_i}}[b_2]$ for
each $i$.

If $b_2<c_i\overline{s_i}^Xd_i$, then
$[c_is_id_i]_{\overline{s_i}}[b_2]$ is quasi-normal $s$-word with
$\beta_ic_i\overline{s_i}d_ib_2\prec w$.

If $b_2>c_i\overline{s_i}^Xd_i$, then
$[c_is_id_i]_{\overline{s_i}}[b_2]=-[b_2][c_is_id_i]_{\overline{s_i}}$
and $[b_2][c_is_id_i]_{\overline{s_i}}$ is quasi-normal $s$-word
with $\beta_ib_2c_i\overline{s_i}d_i\prec
\beta_ib_2a\overline{s}b_1\prec \beta_ia\overline{s}b_1b_2=w$.

If $b_2=c_i\overline{s_i}^Xd_i$, then as above, by Lemma \ref{l40}
and induction on $w$ or by assumption,
$\beta_i[c_is_id_i]_{\overline{s_i}}[b_2]\equiv0 \ \ mod(S,w)$.

These show that for each $i$,
$\beta_i[c_is_id_i]_{\overline{s_i}}[b_2]\equiv0\ \ mod(S,w)$.

Therefore, we may assume that $\lfloor
u\rfloor_s=(asb)=[asb_1]_{\bar{s}}[b_2]$, $a,b\in X^*$, where
$b=b_1b_2$ and $a\bar{s}^Xb_1> b_2$.

Noting that for $[asb_1]_{\bar{s}}=s$ or $[asb_1]_{\bar{s}}=[a_1]
[a_2sb_1]_{\bar{s}}$ with $a_2\bar{s}^Xb_1\leq b_2$ or
$[asb_1]_{\bar{s}}= [asb_{11}]_{\bar{s}}[b_{12}]$ with $b_{12}\leq
b_2$, $\lfloor u\rfloor_s$ is already normal. Now we consider the
remained cases.

Case 2.1. Let $[asb_1]_{\bar{s}}=[a_1][a_2sb_1]_{\bar{s}}$ with
$a_1>a_1a_2\bar{s}^Xb_1>a_2\bar{s}^Xb_1> b_2$. Then we have
$$
\lfloor u\rfloor_s=[[a_{1}][a_2sb_1]_{\bar{s}}][b_2]=[[a_{1}][b_2]]
[a_2sb_1]_{\bar{s}}+[a_{1}][[a_2sb_1]_{\bar{s}}[b_{2}]].
$$

We consider the term $[[a_{1}][b_2]][a_2sb_1]_{\bar{s}}$.

By noting that $a_1>b_2$, we may assume that
$[a_1][b_2]=\sum_{u_i\preceq a_1b_2} \alpha_i[u_i]$ where
 $\alpha_i\in {\bf k}, \ u_i\in ALSW(X)$.
We will prove that $[u_i][a_2sb_1]_{\bar{s}}\equiv0 \ \ mod(S,w). $

If $u_i>a_2\overline{s}^Xb_1$, then $[u_i][a_2sb_1]_{\bar{s}}$ is
quasi-normal $s$-word with $w'=u_ia_2\overline{s}b_1\preceq
a_1b_2a_2\overline{s}b_1\prec w=a_1a_2\overline{s}b_1b_2$.

If $u_i<a_2\overline{s}^Xb_1$, then
$[u_i][a_2sb_1]_{\bar{s}}=-[a_2sb_1]_{\bar{s}}[u_i]$ and
$[a_2sb_1]_{\bar{s}}[u_i]$ is quasi-normal $s$-word with
$w'=a_2\overline{s}b_1u_i\preceq a_2\overline{s}b_1a_1b_2\prec w$,
since $a_1a_2\overline{s}b_1$ is an ALSW.

If $u_i=a_2\overline{s}^Xb_1$, then as above, by Lemma \ref{l40} and
induction on $w$ or by assumption, $[u_i][a_2sb_1]_{\bar{s}}\equiv0 \ \
mod(S,w)$.

This shows that
$$
\lfloor
u\rfloor_s\equiv[a_{1}][[a_2sb_1]_{\bar{s}}[b_{2}]]
\ \ mod(S,w).
$$

By noting that $ a_1>a_2\bar{s}^Xb_1> b_2$,  the result now follows
from the Case 1.

Case 2.2.  Let $[asb_1]_{\bar{s}}=[asb_{11}]_{\bar{s}}[b_{12}]$ with
$a\bar{s}^Xb_{11}>a\bar{s}^Xb_{11}b_{12}>b_{12}> b_2$. Then we have
\begin{eqnarray*}
\lfloor u\rfloor_s=[[asb_{11}]_{\bar{s}}[b_{12}]][b_2]=
[[asb_{11}]_{\bar{s}}[b_2]][b_{12}]+[asb_{11}]_{\bar{s}}[[b_{12}][b_2]].
\end{eqnarray*}

Let us first deal with $[[asb_{11}]_{\bar{s}}[b_2]][b_{12}]$. Since
 $a\bar{s}b_{11}b_2<a\bar{s}b_{11}b_{12}$, we may apply
induction on $w$ and have that
$$
[[asb_{11}]_{\bar{s}}[b_2]][b_{12}]=\sum\alpha_i\beta_i[c_is_id_i]_{\overline{s_i}}[b_{12}],
$$
where $\beta_ic_i\overline{s_i}d_i\preceq a\overline{s}b_{11}b_2, \ w=a\overline{s}b_{11}b_{12}b_2.$

If $b_{12}<c_i\overline{s_i}^Xd_i$, then
$[c_is_id_i]_{\overline{s_i}}[b_{12}]$ is quasi-normal $s$-word with
$w'=\beta_ic_i\overline{s_i}d_ib_{12}\prec w$.

If $b_{12}>c_i\overline{s_i}^Xd_i$, then
$[c_is_id_i]_{\overline{s_i}}[b_{12}]=-[b_{12}][c_is_id_i]_{\overline{s_i}}$
and $[b_{12}][c_is_id_i]_{\overline{s_i}}$ is a quasi-normal
$s$-word with $w'=\beta_ib_{12}c_i\overline{s_i}d_i\preceq
\beta_ib_{12}a\overline{s}b_{11}b_2\prec
a\overline{s}b_{11}b_{12}b_2=w$.

If $b_{12}=c_i\overline{s_i}^Xd_i$, then as above, by Lemma
\ref{l40} and induction on $w$ or by assumption,
$\beta_i[c_is_id_i]_{\overline{s_i}}[b_{12}]\equiv0 \ \ mod(S,w)$.

These show that
$$
\lfloor u\rfloor_s\equiv
[asb_{11}]_{\bar{s}}[[b_{12}][b_2]] \ \ mod(S,w).
$$

Let $[b_{12}][b_2]=[b_{12}b_2]+\sum_{u_i\prec a_1b_2} \alpha_i[u_i]$
where $\alpha_i\in {\bf k}, \ u_i\in ALSW(X)$. By noting that
$a\overline{s}^Xb_{11}>b_{12}b_2$,  we have
$[asb_{11}]_{\bar{s}}[u_i]\equiv 0\ \ mod(S,w)$ for any $i$.
Therefore,
$$
\lfloor u\rfloor_s\equiv
[asb_{11}]_{\bar{s}}[b_{12}b_2] \ \ mod(S,w).
$$

Noting that $[asb_{11}]_{\bar{s}}[b_{12}b_2]$ is quasi-normal and
now
$(w,\delta_{[asb_{11}]_{\bar{s}}[b_{12}b_2]})<(w,\delta_{[asb_{1}]_{\bar{s}}[b_2]})$,
the result follows  by induction.

The proof is complete. $\Box$

\begin{lemma}\label{6} Let $S$ be a ${\bf k}$-monic subset of $Lie_{{\bf k}[Y]}(X)$
in which each multiplication composition is trivial. Then any
element of the ${\bf k}[Y]$-ideal generated by $S$ can be written as
a ${\bf k}[Y]$-linear combination of normal $S$-words.
\end{lemma}

{\bf Proof.} Note that for any $h\in Id(S)$, $h$ can be presented by
a ${\bf k}[Y]$-linear combination of $S$-words of the form
\begin{equation}\label{e6}
(u)_{s}=[c_1][c_2]\cdots[c_{k}]s[d_1][d_2]\cdots[d_{l}]
\end{equation}
with some placement of parentheses, where $s\in S$, $c_j,d_j\in
ALSW(X), k,l\geq 0$. By Lemma \ref{l4} it suffices to prove that
(\ref{e6}) is a linear combination of quasi-normal $S$-words. We
will prove the result by induction on $k+l$. It is trivial when
$k+l=0$, i.e., $(u)_{s}=s$. Suppose that the result holds for
$k+l=n$. Now let us consider
$$
(u)_{s}=[c_{n+1}]([c_1][c_2]\cdots[c_{k}]s[d_1][d_2]\cdots[d_{n-k}])=[c_{n+1}](v)_s.
$$
By inductive hypothesis, we may assume without loss of generality
that $(v)_{s}$ is a quasi-normal $s$-word, i.e., $(v)_{s}=\lfloor
v\rfloor_s=(csd)$ where $c\bar{s}d\in T_A, c,d\in X^*$. If
$c_{n+1}>c\bar{s}^Xd$, then  $(u)_{s}$ is quasi-normal. If
$c_{n+1}<c\bar{s}^Xd$ then $(u)_{s}=-\lfloor v\rfloor_s[c_{n+1}]$
where $\lfloor v\rfloor_s[c_{n+1}]$ is quasi-normal. If
$c_{n+1}=c\bar{s}^Xd$ then by Lemma \ref{l4},
$(u)_{s}=[c_{n+1}](csd)\equiv [c_{n+1}][csd]_{\bar{s}}$. Now the
result follows from the multiplication composition and Lemma
\ref{l40}. $\Box$

\begin{lemma}\label{l8} Let $S$ be a ${\bf k}$-monic subset of $Lie_{{\bf k}[Y]}(X)$
in which each multiplication composition is trivial. Then for any
quasi-normal $S$-word $\lfloor
asb\rfloor_{s}=[a_1][a_2]\cdots[a_k]\lfloor
v\rfloor_s[b_1][b_2]\cdots[b_l]$ with some placement of parentheses,
the  three following $S$-words are linear combinations of normal
$S$-words with the leading words less than $a\bar{s}b$:
\begin{enumerate}
\item[(i)]
$w_1=\lfloor asb\rfloor_{s}|_{[a_i]\mapsto[c]}$ where $c\prec a_i$;
\item[(ii)] $w_2=\lfloor asb\rfloor_{s}|_{[b_j]\mapsto[d]}$ where $d\prec b_j$;
\item[(iii)]$w_3=\lfloor asb\rfloor_{s}|_{\lfloor v\rfloor_s\mapsto\lfloor v'\rfloor_s}$
 where $\overline{\lfloor v'\rfloor}_s\prec\overline{\lfloor
 v\rfloor}_s$.
\end{enumerate}
\end{lemma}
{\bf Proof.} We first prove $(iii)$. For $k+l=1$, for example,
$\lfloor asb\rfloor_{s}=\lfloor v\rfloor_s[b_1]$, it is easy
to see that the result follows from Lemmas \ref{6} and \ref{l40}
since either $\lfloor v'\rfloor_s[b_1]$ or $[b_1]\lfloor
v'\rfloor_s$ is quasi-normal or $w_3$ is the multiplication
composition. Now the result follows by induction on $k+l$.

We now prove $(i)$, and $(ii)$ is similar. For $k+l=1$, $\lfloor
asb\rfloor_{s}=[a_1]\lfloor v\rfloor_s$ and then $w_1=[c]\lfloor
v\rfloor_s$. Then either $\lfloor v\rfloor_s[c]$ or $[c]\lfloor
v\rfloor_s$ is quasi-normal or $w_1$ is equivalent to the
multiplication composition with respect to $w=\overline{\lfloor
v\rfloor}_s^X\overline{\lfloor v\rfloor}_s$. Again by Lemmas \ref{6}
and \ref{l40}, the result holds. For $k+l\geq2$, it follows from
(iii). $\Box$

 \ \

Let $s_1,s_2\in Lie_{{\bf k}[Y]}(X)$ be two ${\bf k}$-monic
polynomials in $Lie_{{\bf k}[Y]}(X)$. If
$a\bar{s}_1^Xb\bar{s}_2^Xc\in ALSW(X)$ for some $a,b,c\in X^*$, then
by Lemma \ref{chi'lem}, there exits a bracketing way
$[a\bar{s}_1^Xb\bar{s}_2^Xc]_{\bar{s}_1^X,\bar{s}_2^X}$ such that
$\overline{[a\bar{s}_1^Xb\bar{s}_2^Xc]}_{\bar{s}_1^X,\bar{s}_2^X}=a\bar{s}_1^Xb\bar{s}_2^Xc$.
Denote
\begin{eqnarray*}
&&[as_1b\bar{s}_2c]_{\bar{s}_1,\bar{s}_2}=\bar{s}_2^Y[a\bar{s}_1^Xb\bar{s}_2^Xc]_{\bar{s}_1^X,\bar{s}_2^X}|_{[\bar{s}_1^X]\mapsto
s_1},\\
&&[a\bar{s}_1bs_2c]_{\bar{s}_1,\bar{s}_2}=\bar{s}_1^Y[a\bar{s}_1^Xb\bar{s}_2^Xc]_{\bar{s}_1^X,\bar{s}_2^X}|_{[\bar{s}_2^X]\mapsto
s_2},\\
&&[as_1bs_2c]_{\bar{s}_1,\bar{s}_2}=[a\bar{s}_1^Xb\bar{s}_2^Xc]_{\bar{s}_1^X,\bar{s}_2^X}|_{[\bar{s}_1^X]\mapsto
s_1,[\bar{s}_2^X]\mapsto s_2}.
\end{eqnarray*}
Thus, the leading words of the above three polynomials are
$a\bar{s}_1b\bar{s}_2c=\bar{s}_1^Y\bar{s}_2^Ya\bar{s}_1^Xb\bar{s}_2^Xc$.

The following lemma is also essential in this paper.

\begin{lemma}\label{l5}Let $S$ be a Gr\"{o}bner-Shirshov basis in
$Lie_{{\bf k}[Y]}(X)$. For any $s_1,s_2\in S, \
\beta_1,\beta_2\in[Y], a_1,a_2,b_1,b_2\in X^*$ such that
$w=\beta_1a_1\bar{s}_1b_1=\beta_2a_2\bar{s}_2b_2\in T_A$, we have
$$
\beta_1[a_1s_1b_1]_{\bar{s}_1}\equiv\beta_2[a_2s_2b_2]_{\bar{s}_2} \ mod(S,w).
$$
\end{lemma}

{\bf Proof.} Let $L$ be  the least common multiple of $\bar{s}_1^Y$
and $\bar{s}_2^Y$. Then
$w^Y=\beta_1\bar{s}_1^Y=\beta_2\bar{s}_2^Y=Lt$ for some $t\in
[Y]$, $w^X=a_1\bar{s}_1^Xb_1=a_2\bar{s}_2^Xb_2$ and
$$
\beta_1[a_1s_1b_1]_{\bar{s}_1}-\beta_2[a_2s_2b_2]_{\bar{s}_2}
=t(\frac{L}{\bar{s}_1^Y}[a_1s_1b_1]_{\bar{s}_1}-\frac{L}{\bar{s}_2^Y}[a_2s_2b_2]_{\bar{s}_2}).
$$

Consider the first case in which $\bar{s}_2^X$ is a subword of $b_1$, i.e.,
$w^X=a_1\bar{s}_1^Xa\bar{s}_2^Xb_2$ for some $a\in X^*$ such that
$b_1=a\bar{s}_2^Xb_2$ and $a_2=a_1\bar{s}_1^Xa$. Then
\begin{eqnarray*}
&&\beta_1[a_1s_1b_1]_{\bar{s}_1}-\beta_2[a_2s_2b_2]_{\bar{s}_2}\\
&=&t(\frac{L}{\bar{s}_1^Y}[a_1s_1a\bar{s}_2^Xb_2]_{\bar{s}_1}-
\frac{L}{\bar{s}_2^Y}[a_1\bar{s}_1^Xas_2b_2]_{\bar{s}_2})\\
&=&tC_3\langle s_1, s_2\rangle_{w'}
\end{eqnarray*}
if $L\neq\bar{s}_1^Y\bar{s}_2^Y$, where $w'=Lw^X$. Since $S$ is a
Gr\"{o}bner-Shirshov basis, $C_3\langle s_1, s_2\rangle\equiv0 \ \
mod(S,Lw^X)$. The result follows from $w=tLw^X=tw'$.

Suppose that $L=\bar{s}_1^Y\bar{s}_2^Y$. By noting that
$\frac{1}{\bar{s}_1^Y}[a_1\bar{s}_1as_2b_2]_{\bar{s}_1,\bar{s}_2}$
and
$\frac{1}{\bar{s}_2^Y}[a_1{s}_1a\bar{s}_2b_2]_{\bar{s}_1,\bar{s}_2}$
are quasi-normal, by Lemma \ref{l4} we have
\begin{eqnarray*}
&&[a_1s_1a\bar{s}_2b_2]_{\bar{s}_1,\bar{s}_2}\equiv\bar{s}_2^Y[a_1s_1a\bar{s}_2^Xb_2]_{\bar{s}_1} \ \ mod(S,w'),\\
&&[a_1\bar{s}_1as_2b_2]_{\bar{s}_1,\bar{s}_2}\equiv\bar{s}_1^Y[a_1\bar{s}^X_1as_2b_2]_{\bar{s}_2}
\ \ mod(S,w').
\end{eqnarray*}
Thus, by Lemma \ref{l8}, we have
\begin{eqnarray*}
&&\beta_1[a_1s_1b_1]_{\bar{s}_1}-\beta_2[a_2s_2b_2]_{\bar{s}_2}\\
&=&t(\bar{s}_2^Y[a_1s_1a\bar{s}_2^Xb_2]_{\bar{s}_1}-
\bar{s}_1^Y[a_1\bar{s}_1^Xas_2b_2]_{\bar{s}_2})\\
&=&t((\bar{s}_2^Y[a_1s_1a\bar{s}_2^Xb_2]_{\bar{s}_1}-[a_1s_1a\bar{s}_2b_2]_{\bar{s}_1,\bar{s}_2})
+([a_1s_1as_2b_2]_{\bar{s}_1,\bar{s}_2}-[a_1s_1a\bar{s}_2b_2]_{\bar{s}_1,\bar{s}_2})\\
&&-([a_1s_1as_2b_2]_{\bar{s}_1,\bar{s}_2}-[a_1\bar{s}_1as_2b_2]_{\bar{s}_1,\bar{s}_2})
-(\bar{s}_1^Y[a_1\bar{s}_1^Xas_2b_2]_{\bar{s}_2}-[a_1\bar{s}_1as_2b_2]_{\bar{s}_1,\bar{s}_2}))\\
&=&t((\bar{s}_1^Y[a_1s_1a\bar{s}_2^Xb_2]_{\bar{s}_1}-[a_1s_1a\bar{s}_2b_2]_{\bar{s}_1,\bar{s}_2})
+[a_1(s_1-[\bar{s}_1])as_2b_2]_{\bar{s}_1,\bar{s}_2}\\
&&-[a_1s_1a(s_2-[\bar{s}_2])b_2]_{\bar{s}_1,\bar{s}_2}
-(\bar{s}_1^Y[a_1\bar{s}_1^Xas_2b_2]_{\bar{s}_2}-[a_1\bar{s}_1as_2b_2]_{\bar{s}_1,\bar{s}_2}))\\
&\equiv& 0 \ \ mod(S,w).
\end{eqnarray*}

 Second, if $\bar{s}_2^X$ is a subword of $\bar{s}_1^X$, i.e.,
$\bar{s}_1^X=a\bar{s}_2^Xb$ for some $a,b\in X^*$, then
$[a_2s_2b_2]_{\bar{s}_2}=[a_1as_2bb_1]_{\bar{s}_2}$. Let
$w'=L\bar{s}_1^X$. Thus, by noting that
$[a_1[as_2b]_{\bar{s}_2}b_1]$ is quasi-normal and by Lemmas \ref{l4}
and \ref{l8},
\begin{eqnarray*}
&&\beta_1[a_1s_1b_1]_{\bar{s}_1}-\beta_2[a_2s_2b_2]_{\bar{s}_2}\\
&=&t(\frac{L}{\bar{s}_1^Y}[a_1s_1b_1]_{\bar{s}_1}-
\frac{L}{\bar{s}_2^Y}[a_1as_2bb_1]_{\bar{s}_2}) \\
&=&t(\frac{L}{\bar{s}_1^Y}[a_1s_1b_1]_{\bar{s}_1}-
\frac{L}{\bar{s}_2^Y}[a_1s_1b_1]_{\bar{s}_1}|_{s_1\mapsto
[as_2b]_{\bar{s}_2}})
-\frac{L}{\bar{s}_2^Y}([a_1as_2bb_1]_{\bar{s}_2}-
[a_1s_1b_1]_{\bar{s}_1}|_{s_1\mapsto[as_2b]_{\bar{s}_2}}) \\
&=&t[a_1(\frac{L}{\bar{s}_1^Y}s_1-\frac{L}{\bar{s}_2^Y}[as_2b]_{\bar{s}_2})b_1]-
\frac{L}{\bar{s}_2^Y}([a_1as_2bb_1]_{\bar{s}_2}-[a_1^X[as_2b]_{\bar{s}_2}b_1])\\
&=&t[a_1C_1\langle s_1,s_2\rangle_{w'} b_1]-\frac{L}{\bar{s}_2^Y}
([a_1as_2bb_1]_{\bar{s}_2}-[a_1[as_2b]_{\bar{s}_2}b_1])\\
&\equiv&0 \ mod(S,w).
\end{eqnarray*}

One more case is possible: A proper suffix of $\bar{s}_1^X$ is a
proper prefix of $\bar{s}_2^X$, i.e., $\bar{s}_1^X=ab$ and
$\bar{s}_2^X=bc$ for some $a,b,c\in X^*$ and $b\neq 1$. Then  $abc$
is an ALSW. Let $w'=Labc$. Then by Lemmas \ref{l4} and \ref{l8}, we
have
\begin{eqnarray*}
&&\beta_1[a_1s_1b_1]_{\bar{s}_1}-\beta_2[a_2s_2b_2]_{\bar{s}_2}\\
&=&t(\frac{L}{\bar{s}_1^Y}[a_1s_1cb_2]_{\bar{s}_1}-\frac{L}{\bar{s}_2^Y}[a_1as_2b_2]_{\bar{s}_2})\\
&=&t\frac{L}{\bar{s}_1^Y}([a_1s_1cb_2]_{\bar{s}_1}-[a_1[s_1c]_{\bar{s}_1}b_2])
   -t\frac{L}{\bar{s}_2^Y}([a_1as_2b_2]_{\bar{s}_2}-[a_1[as_2]_{\bar{s}_2}b_2])\\
   &&+t([a_1C_2\langle s_1,s_2\rangle_{w'}b_2]\\
&\equiv&0 \ mod(S,w).
\end{eqnarray*}

The proof is complete. \ \ \ \ $\Box$

\begin{theorem}{\bf (Composition-Diamond lemma for $Lie_{{\bf k}[Y]}(X)$) } \label{cdL2} Let
$S\subset{Lie_{{\bf k}[Y]}(X)}$ be a nonempty set of ${\bf k}$-monic
polynomials and $Id(S)$ be the ${\bf k}[Y]$-ideal of $Lie_{{\bf
k}[Y]}(X)$ generated by $S$. Then the following statements are
equivalent.
\begin{enumerate}
\item[(i)] $S$ is a Gr\"{o}bner-Shirshov basis in
$Lie_{{\bf k}[Y]}(X)$.
\item[(ii)] $f\in{Id(S)}\Rightarrow{\bar{f}=\beta a\bar{s}b\in T_A}$ for
some $s\in{S},  \ \beta\in[Y]$ and $a,b\in{X^*}$.
\item[(iii)]$Irr(S)=\{[u] \ | \ [u]\in T_N, \ u\neq{\beta a\bar{s}b},  \mbox{ for any }
s\in{S},\ \beta\in[Y], \ a,b\in{X^*}\}$ is a $\bf k$-basis for
$Lie_{{\bf k}[Y]}(X|S)=Lie_{{\bf k}[Y]}(X)/Id(S)$.
\end{enumerate}
\end{theorem}

{\bf Proof.} $(i)\Rightarrow (ii)$. Let $S$ be a
Gr\"{o}bner-Shirshov basis and $0\neq f\in Id(S).$ Then
 by Lemma \ref{6} $f$ has an expression
 $f=\sum\alpha_i\beta_i[a_is_ib_i]_{\bar{s_i}}$, where $\alpha_i\in{\bf k}, \
\beta_i\in[Y], \  a_i,b_i\in X^*,  \ s_i\in S$. Denote
 $w_i=\overline{\beta_i[a_is_ib_i]}_{\bar{s_i}},\ i=1,2,\dots$. Then $w_i=\beta_ia_i\bar{s_i}b_i$.
 We may assume without loss of generality that
$$w_1=w_2=\cdots =w_l\succ w_{l+1}\succeq w_{l+2}\succeq \cdots$$
\noindent for some $l\geq 1$.

The claim of the theorem is obvious if $l=1$.

Now suppose that $l>1$. Then
$\beta_1a_1\bar{s_1}b_1=w_1=w_2=\beta_2a_2\bar{s_2}b_2$. By Lemma \ref{l5},
\begin{eqnarray*}
&&\alpha_1\beta_1[a_1s_1b_1]_{\bar{s_1}}+\alpha_2\beta_2[a_2s_2b_2]_{\bar{s_2}}\\
&=&(\alpha_1+\alpha_2)\beta_1[a_1s_1b_1]_{\bar{s_1}}+\alpha_2(\beta_2[a_2s_2b_2]_{\bar{s_2}}-\beta_1[a_1s_1b_1]_{\bar{s_1}})\\
&\equiv&(\alpha_1+\alpha_2)\beta_1[a_1s_1b_1]_{\bar{s_1}} \ \ \ \
mod(S,w_1).
\end{eqnarray*}

Therefore, if $\alpha_1+\alpha_2\neq 0$ or $l>2$, then the result
follows from the induction on $l$. For the case $\alpha_1+\alpha_2=
0$ and $l=2$, we use the induction on $w_1$. Now the result follows.

$(ii)\Rightarrow (iii).$ For any $f\in Lie_{{\bf k}[Y]}(X)$, we have
\begin{equation*}
f=\sum\limits_{\overline{\beta_i[a_is_ib_i]}_{\bar{s_i}}\preceq \bar
f}\alpha_i\beta_i[a_is_ib_i]_{\bar{s_i}}+
\sum\limits_{\overline{[u_j]}\preceq \bar f}\alpha'_j[u_j],
\end{equation*}
where  $\alpha_i,\alpha'_j\in {\bf k}, \ \beta_i\in[Y], \  \ [u_j]
\in Irr(S)$ and ${s_i}\in S$. Therefore, the set $Irr(S)$ generates
the algebra $Lie_{{\bf k}[Y]}(X)/Id(S)$.

On the other hand, suppose that $h=\sum\alpha_i[u_i]=0$ in
$Lie_{{\bf k}[Y]}(X)/Id(S)$, where $\alpha_i\in {\bf k}$, $[u_i]\in
{Irr(S)}$. This means that $h\in{Id(S)}$. Then all $\alpha_i$ must
be equal to zero. Otherwise, $\overline{h}=u_j$ for some $j$ which
contradicts (ii).
 \\

$(iii)\Rightarrow (i).$
 For any $f,g\in{S}$,  we
have
$$
C_\tau(f,g)_{w}=\sum\limits_{\overline{\beta_i[a_is_ib_i]}_{\bar{s_i}}\prec
w}\alpha_i\beta_i[a_is_ib_i]_{\bar{s_i}}+
\sum\limits_{\overline{[u_j]}\prec w}\alpha'_j[u_j].
$$
For $\tau=1,2,3,4$, since $C_\tau(f,g)_{w}\in {Id(S)}$ and by
$(iii)$, we have
$$
C_\tau(f,g)_{w}=\sum\limits_{\overline{\beta_i[a_is_ib_i]}_{\bar{s_i}}\prec w}\alpha_i\beta_i[a_is_ib_i]_{\bar{s_i}}.
$$
Therefore, $S$ is a Gr\"{o}bner-Shirshov basis. \ \ \ \ $\Box$

\section{Applications}
In this section, all algebras (Lie or associative) are understood to
be taken over an associative and commutative ${\bf k}$-algebra $K$
with identity and all associative algebras are assumed to have
identity.

Let  $\mathcal{L}$ be an arbitrary Lie $K$-algebra which is
presented by generators $X$ and  defining relations $S$,\ $
 \mathcal{L}=Lie_{K}(X|S)$. Let
$K$ have a presentation by generators $Y$ and defining relations
$R$,\ $ K={\bf k}[Y|R] $. Let $\succ_Y$ and $\succ_X$ be deg-lex
orderings on $[Y]$ and $X^*$ respectively. Let $RX=\{rx|r\in R,x\in
X\}$. Then as ${\bf k}[Y]$-algebras,
$$
 \mathcal{L}=Lie_{{\bf k} [Y|R]}(X|S)\cong Lie_{{\bf k}[Y]}(X|S, RX).
$$

As we know, the Poincare-Birkhoff-Witt theorem cannot be generalized
to Lie algebras over an arbitrary ring (see, for example,
\cite{Grivel}). This implies that not any Lie algebra over a
commutative algebra has a faithful representation in an associative
algebra over the same commutative algebra. Following P.M. Cohn (see
\cite{Grivel}), a Lie algebra with the PBW property is said to be
``special".  The first non-special example was given by A.I.
Shirshov in \cite{Shir53} (see also \cite{Shir3}), and he also
suggested that if no nonzero element of $K$ annihilates an absolute
zero-divisor, then a faithful representation always exits. Another
classical non-special example was given  by  P. Cartier
\cite{Cartier}. In the same paper, he proved that each Lie algebra
over Dedekind domain is special. In both examples the Lie algebras
are taken over commutative algebras over $GF(2)$. Shirshov and
Cartier  used ad hoc methods to prove that some elements of
corresponding Lie algebras are not zero though they are zero in the
universal enveloping algebras. P.M. Cohn \cite{Conh} proved that any
Lie algebra over $_{\bf k}K$, where $char ({\bf k})=0$, is special.
Also he claimed that he gave an example of non-special Lie algebra
over a truncated
 polynomial algebra over a filed of characteristic $p>0$.
But he did not give a proof.

Here we find Gr\"{o}bner-Shirshov bases of Shirshov and Cartier's
Lie algebras and then use Theorem \ref{cdL2} to get the results and
we give proof for P.M. Cohn's example of characteristics $2, \ 3$
and $5$. We present an algorithm that one can check for any $p$,
whether Cohn's conjecture is valid.

Note that if $\mathcal{L}=Lie_K(X|S)$, then the universal enveloping
algebra of $\mathcal{L}$ is $U_K(\mathcal{L})=K\langle
X|S^{(-)}\rangle$ where $S^{(-)}$ is just $S$ but substituting all
$[u,v]$ by $uv-vu$.

\begin{example}(Shirshov \cite{Shir53,Shir3})
Let the field ${\bf k}=GF(2)$ and $K={\bf k}[Y|R]$, where
$$Y=\{y_i,i=0,1,2,3\}, \ R=\{y_0y_i=y_i \ (i=0,1,2,3), \ y_iy_j=0 \ (i,j\neq0)\}.$$
Let $\mathcal{L}=Lie_K(X|S_1, S_2)$, where $X=\{x_i, 1\leq i\leq
13\},$ $S_1$ consists of the following relations
\begin{eqnarray*}
&&[x_2,x_1]=x_{11}, \ [x_3,x_1]=x_{13}, \ [x_3,x_2]=x_{12}, \\
&&[x_5,x_3]=[x_6,x_2]=[x_8,x_1]=x_{10}, \\
&& [x_i,x_j]=0 \ \ \ (\mbox{for any other} \ i>j),
\end{eqnarray*}
and $S_2$ consists of the following relations
\begin{eqnarray*}
&& y_0x_i=x_i  \ (i=1,2,\ldots,13),\\
&& x_4=y_1x_1, \ x_5=y_2x_1, \ x_5=y_1x_2, \ x_6=y_3x_1, \
x_6=y_1x_3, \\
&& x_7=y_2x_2, \ x_8 =y_3x_2, \ x_8 =y_2x_3,  \ x_9=y_3x_3, \\
&& y_3x_{11}=x_{10}, \ y_1x_{12}=x_{10}, \ y_2x_{13}=x_{10},\\
&& y_1x_k=0 \ (k=4,5,\ldots,11,13), \ y_2x_t=0 \ (t=4,5,\ldots,12),
\ y_3x_l=0 \ (l=4,5,\ldots,10,12,13).
\end{eqnarray*}
Then $\mathcal{L}$ is not special.
\end{example}

{\bf Proof.} $\mathcal{L}=Lie_K(X|S_1, S_2)=Lie_{{\bf
k}[Y]}(X|S_1,S_2,RX)$. We order $Y$ and $X$ by $y_i>y_j \ \mbox{if}
\ i>j$ and $x_i>x_j \ \mbox{if} \ i>j$ respectively. It is easy to
see that for the ordering  $\succ$ on $[Y]X^*$ as before, $S=S_1\cup
S_2\cup RX\cup\{ y_1x_2=y_2x_1, \ y_1x_3=y_3x_1, \ y_2x_3=y_3x_2\}$
is a Gr\"{o}bner-Shirshov basis in $Lie_{{\bf k}[Y]}(X)$. Since
$x_{10}\in Irr(S)$ and  $Irr(S)$ is a ${\bf k}$-basis of
$\mathcal{L}$ by Theorem \ref{cdL2}, $x_{10}\neq0$ in $\mathcal{L}$.

On the other hand, the universal enveloping algebra of $\mathcal{L}$ has a presentation:
$$
U_K(\mathcal{L})=K\langle X|S_1^{(-)},S_2\rangle\cong{\bf k}[Y]\langle X|S_1^{(-)},S_2,RX\rangle,
$$
where $S_1^{(-)}$ is just $S_1$ but substituting all $[uv]$ by
$uv-vu$.

But the Gr\"{o}bner-Shirshov complement (see Mikhalev-Zolotyhk
\cite{MZ}) of $S_1^{(-)}\cup S_2\cup RX$  in ${\bf k}[Y]\langle
X\rangle$ is
$$S^C=S_1^{(-)}\cup S_2\cup RX\cup\{y_1x_2=y_2x_1, \
y_1x_3=y_3x_1,  \ y_2x_3=y_3x_2, \  x_{10}=0\}.$$

Thus, $\mathcal{L}$ is not special. \ \ \ \ $\Box$

\begin{example} (Cartier  \cite{Cartier}) Let ${\bf k}=GF(2)$, $K={\bf k}[y_1,y_2,y_3|y_i^2=0,\ i=1,2,3]$ and
$\mathcal{L}=Lie_{K}(X|S)$, where $X=\{x_{ij},1\leq i\leq j\leq3\}$ and
$$
S=\{[x_{ii},x_{jj}]=x_{ji} \ (i>j), [x_{ij},x_{kl}]=0 \
(\mbox{otherwise}), \ y_3x_{33}=y_2x_{22}+y_1x_{11}\}.
$$
Then $\mathcal{L}$ is not special.
\end{example}
{\bf Proof.} Let $Y=\{y_1,y_2,y_3\}$. Then
$$
\mathcal{L}=Lie_{K}(X|S)\cong Lie_{{\bf k}[Y]}(X|S,
y_i^2x_{kl}=0\ (\forall i,k,l)).
$$
Let $y_i>y_j \ \mbox{if} \ i>j$ and $x_{ij}>x_{kl} \ \mbox{if} \
(i,j)>_{lex}(k,l)$ respectively. It is easy to see that for the
ordering $\succ$ on $[Y]X^*$ as before, $S'=S\cup \{y_i^2x_{kl}=0\
(\forall i,k,l)\}\cup S_1$ is a Gr\"{o}bner-Shirshov basis in
$Lie_{{\bf k}[Y]}(X)$, where $S_1$ consists of the following
relations
\begin{eqnarray*}
&&y_3x_{23}=y_1x_{12}, \  y_3x_{13}=y_2x_{12}, \ y_2x_{23}=y_1x_{13}, \ y_3y_2x_{22}=y_3y_1x_{11}, \\
&&y_3y_1x_{12}=0, \ y_3y_2x_{12}=0,  \ y_3y_2y_1x_{11}=0,   \ y_2y_1x_{13}=0.
\end{eqnarray*}

The universal enveloping algebra of $\mathcal{L}$ has a
presentation:
$$
U_K(\mathcal{L})=K\langle X|S^{(-)}\rangle\cong{\bf k}[Y]\langle
X|S^{(-)},y_i^2x_{kl}=0\ (\forall i,k,l)\rangle.
$$

In $U_K(\mathcal{L})$, we have (cf. \cite{Cartier})
\begin{eqnarray*}
0=y_3^2x_{33}^2=(y_2x_{22}+y_1x_{11})^2=y_2^2x_{22}^2+y_1^2x_{11}^2+y_2y_1[x_{22},x_{11}]
 = y_2y_1x_{12}.
\end{eqnarray*}
On the other hand, since $y_2y_1x_{12}\in Irr(S')$,
$y_2y_1x_{12}\neq0$ in $\mathcal{L}$. Thus, $\mathcal{L}$ is not
special.  \ \ \ \ $\Box$

\begin{conjecture}(Cohn \cite{Conh}) Let $K={\bf k}[y_1,y_2,y_3|y_i^p=0, i=1,2,3]$ be
the algebra of truncated polynomials over a field ${\bf k}$ of
characteristic $p>0$. Let
$$
\mathcal{L}_p=Lie_{K}(x_1,x_2,x_3 \ | \ y_3x_3=y_2x_2+y_1x_1).
$$
Then $\mathcal{L}_p$ is not special.  We call $\mathcal{L}_p$ the
Cohn's Lie algebra.

\end{conjecture}

\noindent{\bf Remark} (see \cite{Conh}): In $U_K(\mathcal{L}_p)$ we
have
$$
0= (y_3x_3)^p = (y_2x_2)^p + \Lambda_p(y_2x_2, y_1x_1) + (y_1x_1)^p=
\Lambda_p(y_2x_2, y_1x_1),
$$
where $\Lambda_p$ is a Jacobson-Zassenhaus Lie polynomial. P.M. Cohn
conjectured that $\Lambda_p(y_2x_2, y_1x_1)\neq 0 $ in
$\mathcal{L}_p$.

\begin{theorem}
Cohn's Lie  algebras $\mathcal{L}_2$, $\mathcal{L}_3$ and
$\mathcal{L}_5$ are not special.
\end{theorem}
{\bf Proof.} Let $Y=\{y_1,y_2,y_3\}$, $X=\{x_1,x_2,x_3\}$ and
$S=\{y_3x_3=y_2x_2+y_1x_1, \  y_i^px_j=0,\ 1\leq i,j\leq3\}$. Then
$\mathcal{L}_p\cong Lie_{{\bf k}[Y]}(X  | S)$ and
 $U_K(\mathcal{L}_p)\cong{\bf k}[Y]\langle X|S\rangle$. Suppose that $S^C$ is a Gr\"{o}bner-Shirshov complement
 of $S$ in $Lie_{{\bf k}[Y]}(X)$. Let $S_{_{X^{^p}}}\subset\mathcal{L}_p$ be the set of all the elements of
$S^C$ whose $X$-degrees do
 not exceed $p$.

First, we consider $p=2$ and prove the element
$\Lambda_2=[y_2x_2,y_1x_1]=y_2y_1[x_2x_1]\neq0$ in $\mathcal{L}_2$.

Then by Shirshov's algorithm we have that $S_{X^2}$ consists of the
following relations
\begin{eqnarray*}
&&y_3x_3=y_2x_2+y_1x_1, \ y_i^2x_j=0\ (1\leq i,j\leq3), \  y_3y_2x_2=y_3y_1x_1, \ y_3y_2y_1x_1=0, \\
&&y_2[x_3x_2]=y_1[x_3x_1],  \ y_3y_1[x_2x_1]=0, \ y_2y_1[x_3x_1]=0.
\end{eqnarray*}
Thus, $\Lambda_2$ is in the ${\bf k}$-basis $Irr(S^C)$ of
$\mathcal{L}_2$.

Now, by the above remark, $\mathcal{L}_2$ is not special.

Second, we consider $p=3$ and prove the element
$\Lambda_3=y_2^2y_1[x_2x_2x_1]+y_2y_1^2[x_2x_1x_1]\neq0$ in
$\mathcal{L}_3$.

Then again by Shirshov's algorithm, $S_{X^3}$ consists of the
following relations
\begin{eqnarray*}
&&y_3x_3=y_2x_2+y_1x_1, \ y_i^3x_j=0\ (1\leq i,j\leq3), \  y_3^2y_2x_2=y_3^2y_1x_1, \ y_3^2y_2^2y_1x_1=0, \\
&&y_2[x_3x_2]=-y_1[x_3x_1],  \ y_3^2y_1[x_2x_1]=0, \ y_2^2y_1[x_3x_1]=0, \\
&&y_3y_2^2[x_2x_2x_1]=y_3y_2y_1[x_2x_1x_1], \ y_3y_2^2y_1[x_2x_1x_1]=0, \ y_3y_2y_1[x_2x_2x_1]=y_3y_1^2[x_2x_1x_1].
\end{eqnarray*}
Thus, $y_2^2y_1[x_2x_2x_1],y_2y_1^2[x_2x_1x_1]\in Irr(S^C)$, which
implies $\Lambda_3\neq0$ in $\mathcal{L}_3$.

Third, let $p=5$. Again by Shirshov's algorithm, $S_{X^5}$ consists
of the following relations
\begin{eqnarray*}
&&1)\ \ y_3x_3=y_2x_2+y_1x_1, \\
&&2)\ \ y_i^5x_j=0,\ 1\leq i,j\leq3, \\
&&3)\ \ y_3^4y_2x_2=-y_3^4y_1x_1, \\
&&4)\ \  y_3^4y_2^4y_1x_1=0, \\
&&5)\ \ y_2[x_3x_2]=-y_1[x_3x_1], \\
&& 6)\ \ y_3^4y_1[x_2x_1]=0, \\
&&7)\  \  y_2^4y_1[x_3x_1]=0, \\
&&8)\ \  y_3^3y_2^2[x_2x_2x_1]=y_3^3y_2y_1[x_2x_1x_1], \\
&&9)\  \ y_3^3y_2^4y_1[x_2x_1x_1]=0, \\
&&10)\ \  y_3^3y_2y_1[x_2x_2x_1]=y_3^3y_1^2[x_2x_1x_1],\\
&&11)\  \ y_1[x_3x_2x_3x_1]=0, \\
&& 12)\  \ y_1[x_3x_1x_2x_1]=0, \\
&&13)\ \  y_1[x_3x_2x_2x_1]=-y_1[x_3x_2x_1x_2], \\
&&14)\  \  y_2[x_3x_1x_2x_1]=0,\\
&&15)\ \
y_3^2y_2^3[x_2x_2x_2x_1]=2y_3^2y_2^2y_1[x_2x_2x_1x_1]-y_3^2y_2y_1^2[x_2x_1x_1x_1],
\\
&&16)\ \   y_3^3y_2^3y_1^2[x_2x_1x_1x_1]=0,
\\
&&17)\  \  y_3^2y_2^2y_1[x_2x_2x_2x_1]=2y_3^2y_2y_1^2[x_2x_2x_1x_1]-y_3^2y_1^3[x_2x_1x_1x_1],  \\
&&18)\  \  y_3^2y_2^4y_1^2[x_2x_1x_1x_1]=0, \\
\end{eqnarray*}
\begin{eqnarray*}
&&19)\  \  y_3^2y_2^4y_1[x_2x_2x_1x_1]=\frac{1}{2}y_3^2y_2^3y_1^2[x_2x_1x_1x_1], \\
&& 20)\  \  y_3^3y_1^2[x_2x_2x_1x_2x_1]=0, \\
&&21)\  \  y_3^3y_2y_1[x_2x_1x_2x_1x_1]=0, \\
&& 22)\  \  y_3^3y_1^2[x_2x_1x_2x_1x_1]=0,\\
&&23)\  \   y_3^3y_2^2[x_2x_1x_2x_1x_1]=0, \\
&&24)\  \  y_3^2y_2^2y_1[x_2x_2x_1x_2x_1]=-y_3^2y_2y_1^2[x_2x_1x_2x_1x_1], \\
&&25)\  \  y_3^2y_2y_1^2[x_2x_2x_1x_2x_1]=-y_3^2y_1^3[x_2x_1x_2x_1x_1],\\
&&26)\  \  y_3^2y_2^4y_1^2[x_2x_1x_2x_1x_1]=0,\\
&&27)\  \
y_3y_2^4[x_2x_2x_2x_2x_1]=3y_3y_2^3y_1[x_2x_2x_2x_1x_1]-y_3y_2^3y_1[x_2x_2x_1x_2x_1]-
3y_3y_2^2y_1^2[x_2x_2x_1x_1x_1]\\
&&\ \ \ \ \ \ \ \ \  \ \ \ \ \ \ \ \ \ \  \ \ \ \ \ \ \ \ \ \ \  \ \
-2y_3y_2^2y_1^2[x_2x_1x_2x_1x_1]+y_3y_2y_1^3[x_2x_1x_1x_1x_1],
\\
&&28)\  \
y_3y_2^3y_1[x_2x_2x_2x_2x_1]=3y_3y_2^2y_1^2[x_2x_2x_2x_1x_1]-y_3y_2^2y_1^2[x_2x_2x_1x_2x_1]-
3y_3y_2y_1^3[x_2x_2x_1x_1x_1]\\
&&\ \ \ \ \ \ \ \ \  \ \ \ \ \ \ \ \ \ \  \ \ \ \ \ \ \ \ \ \ \  \ \
\ \ -2y_3y_2y_1^3[x_2x_1x_2x_1x_1]+y_3y_1^4[x_2x_1x_1x_1x_1],
\\
&&29)\  \  y_3y_2^4y_1^3[x_2x_1x_1x_1x_1]=0, \\
&&30)\  \  y_3^2y_2^3y_1^3[x_2x_1x_1x_1x_1]=0,\\
&&31)\  \ y_3y_2^4y_1^2[x_2x_2x_1x_1x_1]=-\frac{2}{3}y_3y_2^4y_1^2[x_2x_1x_2x_1x_1]+\frac{1}{3}y_3y_2^3y_1^3[x_2x_1x_1x_1x_1],\\
&&32)\  \ y_3y_2^4y_1[x_2x_2x_2x_1x_1]=\frac{1}{3}y_3y_2^4y_1[x_2x_2x_1x_2x_1]+y_3y_2^3y_1^2[x_2x_2x_1x_1x_1]\\
&&\ \ \ \ \ \ \ \ \  \ \ \ \ \ \ \ \ \ \  \ \ \ \ \ \ \ \ \ \  \ \ \
\ \
+\frac{2}{3}y_3y_2^3y_1^2[x_2x_1x_2x_1x_1]-\frac{1}{3}y_3y_2^2y_1^3[x_2x_1x_1x_1x_1],\\
&&33)\  \ y_2^3y_1^2[x_3x_3x_1x_3x_1]=0, \\
&&34)\  \  y_2^3y_1^2[x_3x_1x_3x_1x_1]=0, \\
&&35)\  \ y_3^3y_2^2y_1^3[x_2x_1x_1x_1x_1]=0,\\
&&36)\  \
y_3^2y_2^3y_1^2[x_2x_2x_1x_1x_1]=-\frac{2}{3}y_3^2y_2^3y_1^2[x_2x_1x_2x_1x_1]+\frac{2}{3}y_3^2y_2^2y_1^3[x_2x_1x_1x_1x_1].
\end{eqnarray*}
Thus, $\overline{\Lambda_5(y_2x_2,
y_1x_1)}=y_2^4y_1[x_2x_2x_2x_2x_1]\in Irr(S^C)$, which implies
$\Lambda_5\neq0$ in $\mathcal{L}_5$.
 \ \ \ \ $\Box$
\\

\noindent{\bf Remarks:} Note that the Jacobson-Zassenhaus Lie
polynomial $\Lambda_p(y_2x_2, y_1x_1)$ is of $X$-degree $p$. Then
$\overline{\Lambda_p(y_2x_2, y_1x_1)}\in Irr(S^C)$ if and only if
$\overline{\Lambda_p(y_2x_2,
y_1x_1)}\in Irr(S_{X^p})$. Since the defining relation of $\mathcal{L}_p$ is homogenous on $X$, $S_{X^p}$ is a finite set. By Shirshov's algorithm, one can compute $S_{X^p}$ for $\mathcal{L}_p$.\\

Now we give some examples which are special Lie algebras.

\begin{lemma}\label{l10} Suppose that $f$ and $g$ are two  polynomials
in $Lie_{{\bf k}[Y]}(X)$ such that $f$ is ${\bf k}[Y]$-monic and
$g=rx$, where $r\in {\bf k}[Y]$ and $x\in X$, is ${\bf k}$-monic.
 Then each inclusion composition of $f$ and $g$ is trivial modulo $\{f\}\cup rX$.
\end{lemma}

{\bf Proof.} Suppose that $\bar{f}=[axb]$ for some $a,b\in X^*$, $f=\bar{f}+f'$
 and $g=\bar{r}x+r'x$. Then $w=\bar{r}axb$ and
\begin{eqnarray*}
C_1\langle f, g\rangle_w
&=&\bar{r}f-[a[rx]b]_{\bar{r}x}\\
&=&\bar{r}f'-r'[axb]\\
&=& rf'-r'f\\
&\equiv&0 \ \ mod(\{f\}\cup rX, w). \ \ \ \Box
\end{eqnarray*}

\begin{theorem}\label{t4.5}  For an arbitrary commutative ${\bf k}$-algebra $K={\bf k}[Y|R]$, if $S$ is
a  Gr\"{o}bner-Shirshov basis in $Lie_{{\bf k}[Y]}(X)$ such that for
any $s\in S$, $s$ is ${\bf k}[Y]$-monic, then
$\mathcal{L}=Lie_{K}(X|S)$ is  special.
\end{theorem}
{\bf Proof.} Assume without loss of generality that  $R$ is a
Gr\"{o}bner-Shirshov basis
 in ${\bf k}[Y]$. Note that $\mathcal{L}\cong Lie_{{\bf k}[Y]}(X|S, RX)$. By Lemma \ref{l10}, $S\cup RX$ is a
Gr\"{o}bner-Shirshov basis in $Lie_{{\bf k}[Y]}(X)$.

On the other hand, in $U_K(\mathcal{L})\cong {\bf k}[Y]\langle
X|S^{(-)}, RX\rangle$, $S^{(-)}\cup RX$ is  a Gr\"{o}bner-Shirshov
basis in ${\bf k}[Y]\langle X\rangle$ in the sense of the paper
\cite{MZ}.

Thus for any $u\in Irr(S\cup RX)$ in $Lie_{{\bf k}[Y]}(X)$,  we have
$\bar{u}\in Irr(S^{(-)}\cup RX)$ in ${\bf k}[Y]\langle X\rangle$.
This completes the proof.  \ \ \ $\Box$

\begin{corollary}\label{co4.6}
Any Lie $K$-algebra $\mathcal{L}=Lie_K(X|f)$ with one monic defining
relation $f=0$ is special.
\end{corollary}
{\bf Proof.} Let $K={\bf k}[Y|R]$, where $R$ is a
Gr\"{o}bner-Shirshov basis
 in ${\bf k}[Y]$. We can regard $f$ as a ${\bf k}[Y]$-monic element in $Lie_{{\bf
k}[Y]}(X)$. Note that any subset of $Lie_{{\bf
k}[Y]}(X)$ consisting of a single ${\bf k}[Y]$-monic element is a
 Gr\"{o}bner-Shirshov basis. Thus by Theorem \ref{t4.5}, $\mathcal{L}=Lie_K(X|f)\cong
 Lie_{{\bf k}[Y]}(X|f, RX)$ is special. \ \ \ $\Box$

\begin{corollary} (\cite{Birkhoff, Witt})
 If $\mathcal{L}$ is a free $K$-module, then $\mathcal{L}$ is special.
\end{corollary}
{\bf Proof.} Let $X=\{x_i, \ i\in I\}$ be a $K$-basis of
$\mathcal{L}$ and $[x_i,x_j]=\sum\alpha_{ij}^lx_l$, where
$\alpha_{ij}^l\in K$ and $i,j\in I$. Then
$\mathcal{L}=Lie_K(X|[x_i,x_j]-\sum\alpha_{ij}^lx_l, \ i>j,\ i,j\in
I)$. Suppose that $K={\bf k}[Y|R]$, where $R$ is a
Gr\"{o}bner-Shirshov basis
 in ${\bf k}[Y]$. Since $S=\{[x_i,x_j]-\sum\alpha_{ij}^lx_l, \ i>j,\ i,j\in I\}$
 is a ${\bf k}[Y]$-monic Gr\"{o}bner-Shirshov basis
  in $Lie_{{\bf
k}[Y]}(X)$,  by Theorem \ref{t4.5}, $\mathcal{L}= Lie_{K}(X|S) \cong
Lie_{{\bf k}[Y]}(X|S, RX)$ is special. \ \ \ $\Box$

Now we give other applications.

\begin{theorem} Suppose that $S$ is a finite homogeneous subset of $Lie_{{\bf k}}(X)$.
Then the word problem of $Lie_{K}(X|S)$ is solvable for any finitely
generated commutative ${\bf k}$-algebra $K$.
\end{theorem}
{\bf Proof.} Let $S^C$ be a Gr\"{o}bner-Shirshov complement of $S$
in $Lie_{{\bf k}}(X)$. Clearly, $S^C$ consists of homogeneous
elements in $Lie_{{\bf k}}(X)$ since the compositions of homogeneous
elements are homogeneous. Since $K$ is finitely generated
commutative ${\bf k}$-algebra, we may assume that $K={\bf k}[Y|R]$
with $R$ a finite Gr\"{o}bner-Shirshov basis in ${\bf k}[Y]$.
 By Lemma \ref{l10}, $S^C\cup RX$ is a Gr\"{o}bner-Shirshov
basis in $Lie_{{\bf k}[Y]}(X)$.
 For a given $f\in Lie_{K}(X)$, it is obvious that after a finite number of steps one can
  write down all the elements of $S^C$ whose $X$-degrees do not exceed the degree of
  $\bar{f}^X$. Denote the set of such elements by $S_{\bar{f}^X}$. Then $S_{\bar{f}^X}$ is
  a finite set. By Theorem \ref{cdL2}, the result follows.  \ \ \ $\Box$

\begin{theorem} Every finitely or countably generated Lie $K$-algebra can be embedded into a two-generated Lie
$K$-algebra, where $K$ is an arbitrary commutative ${\bf
k}$-algebra.
\end{theorem}

{\bf Proof.} Let $K={\bf k}[Y|R]$ and $\mathcal{L}= Lie_{K }(X|S)$
where $X=\{x_i,i\in I\}$ and $I$ is a subset of the set of nature
numbers. Without loss of generality, we may assume that with the
ordering $\succ$ on $[Y]X^*$ as before, $S\cup RX$ is a
Gr\"{o}bner-Shirshov basis in $Lie_{{\bf k}[Y]}(X)$.

Consider the algebra $ \mathcal{L}'=Lie_{{\bf k}[Y]}(X,a,b|S')$
where $S'=S\cup RX \cup R\{a,b\} \cup\{[aab^iab]-x_i,  i\in I\}$.

Clearly, $ \mathcal{L}'$ is a Lie $K$-algebra generated by $a,b$.
Thus, in order to prove the theorem, by using our Theorem
\ref{cdL2}, it suffices to show that with the ordering $\succ$ on
$[Y](X\cup\{a,b\})^{*}$ as before, where $a\succ b\succ x_i,\ x_i\in
X$, $S'$ is a Gr\"{o}bner-Shirshov basis in $Lie_{{\bf
k}[Y]}(X,a,b)$.

It is clear that all the possible compositions of multiplication,
intersection and inclusion are trivial. We only check the external
compositions of some $f\in S$ and $ra\in Ra$: Let
$w=Lu_1\bar{f}^Xu_2au_3$ where $L=L(\bar{f}^Y, \bar{r})$ and $
u_1\bar{f}^Xu_2au_3\in ALSW(X,a,b)$. Then
\begin{eqnarray*}
&&C_3\langle f, ra\rangle_w\\
&=&\frac{L}{\bar{f}_1^Y}[u_1fu_2au_3]_{\bar{f}}-\frac{L}
{\bar{r}}[u_1\bar{f}^Xu_2(ra)u_3]\\
&=&(\frac{L}{\bar{f}_1^Y}[u_1fu_2au_3]_{\bar{f}}-
r\frac{L}{\bar{r}}[u_1\bar{f}^Xu_2au_3]_{\bar{f}^X})-
(\frac{L}{\bar{r}}[u_1\bar{f}^Xu_2(ra)u_3]- r\frac{L}
{\bar{r}}[u_1\bar{f}^Xu_2au_3]_{\bar{f}^X}) \\
&=&([u_1(\frac{L}{\bar{f}_1^Y}f)u_2au_3]_{\bar{f}}-
[u_1(r\frac{L}{\bar{r}}\bar{f}^X)u_2au_3]_{\bar{f}^X})-
r\frac{L}{\bar{r}}([u_1\bar{f}^Xu_2au_3]- [u_1\bar{f}^Xu_2au_3]_{\bar{f}^X}) \\
&\equiv&[u_1C_3\langle f,rx\rangle_{w'}u_2au_3] \ \ mod(S', w)
\end{eqnarray*}
for some $x$ occurring in $\bar{f}^X$ and $w'=L\bar{f}^X$. Since
$S\cup RX$ is a Gr\"{o}bner-Shirshov basis in $Lie_{{\bf k}[Y]}(X)$,
$C_3\langle f,rx\rangle_{w'}\equiv0  \ \ mod(S\cup RX, w')$. Thus by
Lemma \ref{l8}, $[u_1C_3\langle f,rx\rangle_{w'}u_2au_3]\equiv0  \ \
mod(S', w)$. \ \ \ \ $\Box$
\\


\begin{thebibliography}{11}

\bibitem{AL}William W. Adams and Philippe Loustaunau, An introduction to Gr\"{o}bner
bases, Graduate Studies in Mathematics, Vol. 3, American
Mathematical Society (AMS), 1994.


\bibitem{Be78}G.M. Bergman, The diamond lemma for ring theory, {\it Adv.
Math.}, {\bf 29}(1978), 178-218.


\bibitem{Birkhoff} G. Birkhoff, Representability of Lie algebras and Lie groups by matrices,
{\it Ann. of Math}., \textbf{38}(2)(1937), 526-532. (Selected
Papers, Birkhauser 1987, 332-338.)


\bibitem{Bo72}L.A. Bokut, Insolvability of the word problem for Lie algebras, and
subalgebras of finitely presented Lie algebras, {\it Izvestija AN
USSR (mathem.)}, {\bf 36}(6)(1972), 1173-1219.

\bibitem{Bo76}L.A. Bokut, Imbeddings into simple associative
algebras, {\it Algebra i Logika}, {\bf 15}(1976), 117-142.

\bibitem{bc07} L.A. Bokut and Yuqun Chen,  Gr\"{o}bner-Shirshov bases for Lie algebras: after
A.I. Shirshov,  {\it Southeast Asian Bull. Math.},  {\bf 31}(2007),
1057-1076.

\bibitem{BC}L.A. Bokut and Yuqun Chen, Gr\"{o}bner-Shirshov
 bases: Some new results, Proceedings of the Second International Congress
 in Algebra and Combinatorics, World Scientific, 2008, 35-56.

\bibitem{BCC08}L.A. Bokut, Yuqun Chen and Yongshan Chen, Composition-Diamond lemma
for tensor product of free algebras, \emph{Journal of Algebra},
\textbf{323}(2010), 2520-2537.

\bibitem{BCL08}L.A. Bokut,  Yuqun Chen and Cihua Liu,
Gr\"{o}bner-Shirshov bases for dialgebras,\emph{ International Journal of
Algebra and Computation}, \textbf{20}(3)(2010), 391-415.

\bibitem{BCD08}L.A. Bokut,  Yuqun Chen and Xueming Deng,
Gr\"{o}bner-Shirshov bases for Rota-Baxter algebras, \emph{Siberian Math.
J.}, \textbf{51}(6)(2010), 978-988.

\bibitem{BCLi08}L.A. Bokut,  Yuqun Chen and Yu Li,
Gr\"{o}bner-Shirshov bases for Vinberg-Koszul-Gerstenhaber
 right-symmetric algebras,
 \emph{Fundamental and Applied Mathematics}, \textbf{14}(8)(2008), 55-67
 (in Russian). \emph{Journal of Mathematical Sciences}, \textbf{166}(2010), 603-612.

\bibitem{BCM} L.A. Bokut, Yuqun Chen and Qiuhui Mo,
Gr\"{o}bner-Shirshov bases and embeddings of algebras, \emph{International
Journal of Algebra and Computation}, \textbf{20}(2010), 875-900.

\bibitem{BCQ08}L.A. Bokut, Yuqun Chen and Jianjun Qiu,
Gr\"{o}bner-Shirshov bases for associative algebras with multiple
operators and free Rota-Baxter algebras, {\it Journal of Pure and
Applied Algebra}, {\bf 214}(2010), 89-100.

\bibitem{BFK}L.A. Bokut, Y. Fong and W.-F. Ke, Composition-Diamond lemma for
 associative conformal algebras, \emph{Journal of Algebra},  {\bf 272}(2004), 739-774.

\bibitem{BFKK00}L.A. Bokut, Y. Fong, W.-F. Ke and P.S. Kolesnikov, Gr\"obner and
Gr\"obner-Shirshov bases in algebra and conformal algebras, {\it
Fundamental and Applied Mathematics}, {\bf 6}(3)(2000), 669-706.

\bibitem{BK03}L.A. Bokut and  P.S. Kolesnikov, Gr\"obner-Shirshov bases: from their
incipiency to the present, {\it Journal of Mathematical Sciences},
{\bf 116}(1)(2003),  2894-2916.

\bibitem{BK05}L.A. Bokut and P.S. Kolesnikov, Gr\"obner-Shirshov bases, conformal
algebras and pseudo-algebras, {\it Journal of Mathematical
Sciences}, {\bf 131}(5)(2005), 5962-6003.

\bibitem{BKu94}L.A. Bokut and G. Kukin, Algorithmic and Combinatorial algebra, Kluwer Academic Publ., Dordrecht,
1994.



\bibitem{Bu70}B. Buchberger, An algorithmical criteria for the
solvability of algebraic systems of equations, {\it
Aequationes Math.}, {\bf 4}(1970), 374-383.

\bibitem{BuCL}B. Buchberger, G.E. Collins, R. Loos and R. Albrecht, Computer algebra,
symbolic and algebraic computation, Computing Supplementum, Vol.4,
New York: Springer-Verlag, 1982.

\bibitem{BuW}B. Buchberger and Franz Winkler, Gr\"{o}bner bases and
applications, London Mathematical Society Lecture Note Series,
Vol.251, Cambridge: Cambridge University Press, 1998.


\bibitem{Cartier} P. Cartier, Remarques sur le th$\acute{\mbox e}$or$\grave{\mbox e}$me de
Birkhoff-Witt, \emph{Annali della Scuola Norm. Sup. di Pisa s$\acute{\mbox e}$rie} III vol XII(1958), 1-4.


\bibitem{CFL}K.-T. Chen, R. Fox, R. Lyndon, Free differential calculus IV: The quotient group of the
lower central series, \emph{Ann. of Math}., \textbf{68}(1958) 81-95.



\bibitem{CCZ}Yuqun Chen, Yongshan Chen and Chanyan Zhong, Composition-Diamond
lemma for modules,\emph{ Czechoslovak Math. J.},
\textbf{60}(135)(2010), 59-76.

\bibitem{CJZ}Yuqun Chen, Jing Li and Mingjun Zeng, Composition-Diamond lemma
for non-associative algebras over a polynomial algebra,
\emph{Southeast Asian Bull. Math.},  \textbf{34}(2010), 629-638.

\bibitem{Chi}E.S. Chibrikov, On free Lie conformal algebras, \emph{Vestnik Novosibirsk State University},
\textbf{4}(1), 65-83(2004).

\bibitem{Chi2}E.S. Chibrikov, Lyndon-Shirshov words and the intersection of principal ideals in a free Lie algebra, preprint.

\bibitem{Conh} P.M. Cohn, A remark on the Birkhoff-Witt theorem, \emph{Journal London Math. Soc}., \textbf{38}(1963), 197-203.



\bibitem{CLO} David A. Cox, John Little  and Donal O'Shea, Ideals,
varieties and algorithms: An introduction to computational algebraic
geometry and commutative algebra, Undergraduate Texts in
Mathematics, New York: Springer-Verlag, 1992.

\bibitem{Ei} David Eisenbud, Commutative algebra with a view toward algebraic geometry,
Graduate Texts in Math., Vol.150, Berlin and New York:
Springer-Verlag, 1995.


\bibitem{Grivel}P.-P. Grivel, Une histoire du th$\acute{\mbox e}$or$\grave{\mbox e}$me de
Poincar$\acute{\mbox e}$-Birkhoff-Witt, {\it Expositiones
Mathematicae}, \textbf{22}(2004), 145-184.



\bibitem{Higgins} Philip J. Higgins, Baer Invariants and the Birkhoff-Witt Theorem, \emph{Journal of Algebra}, \textbf{11}(1969),
469-482.



\bibitem{Hi64}H. Hironaka, Resolution of singularities of an algebraic variety
over a field if characteristic zero, I, II, {\it Ann. of Math.},
\textbf{79}(1964), 109-203, 205-326.

\bibitem{KL}S.-J. Kang and K.-H. Lee, Gr\"obner-Shirshov bases for irreducible
$sl_{n+1}$-modules, {\it Journal of Algebra}, \textbf{232}(2000),
1-20.


\bibitem{Lazard52} M. Lazard, Sur les alg$\grave{\mbox e}$bres enveloppantes universelles
de certaines alg$\grave{\mbox e}$bres de Lie, \emph{CRAS Paris},
\textbf{234}(1)(1952), 788-791.


\bibitem{Lazard54} M. Lazard,  Sur les alg$\grave{\mbox e}$bres enveloppantes de certaines
alg$\grave{\mbox e}$bres de Lie, \emph{Publ. Sci. Univ. Alger
s$\acute{\mbox e}$rie A}, \textbf{1}(1954), 281-294.

\bibitem{Lyndon} R.C. Lyndon, On Burnside's problem I, \emph{Trans. Amer. Math. Soc.},  \textbf{77}(1954), 202-215.


\bibitem{Mik89}A.A. Mikhalev, The junction lemma and the equality
problem for color Lie superalgebras, {\it Vestnik Moskov. Univ. Ser.
I Mat. Mekh.}, {\bf 5}(1989), 88-91. English translation: {\it
Moscow Univ. Math. Bull.}, {\bf 44}(1989), 87-90.

\bibitem{Mik92}A.A. Mikhalev, The composition lemma for color Lie
superalgebras and for Lie $p$-superalgebras, {\it Contemp. Math.},
{\bf 131}(2)(1992), 91-104.

\bibitem{Mik96}A.A. Mikhalev, Shirshov's composition techniques in
Lie superalgebra (non-commutative Gr\"{o}bner bases). {\it Trudy.
Sem. Petrovsk.}, {\bf 18}(1995), 277-289. English translation: {\it
J. Math. Sci.}, {\bf 80}(1996), 2153-2160.

\bibitem{MZ}A.A. Mikhalev and A.A. Zolotykh, Standard Gr\"obner-Shirshov bases of free algebras over
rings,
I. Free associative algebras,  \emph{International Journal of
Algebra and Computation}, \textbf{8}(6)(1998), 689-726.

\bibitem{NouazeRevoy} Yvon Nouaze and  Philippe Revoy,
Un eas particulier du th$\acute{\mbox e}$or$\grave{\mbox e}$me de Poincar$\acute{\mbox e}$-Birkhoff-Witt.
\emph{CRAS Paris}, s$\acute{\mbox e}$rie A(1971), 329-331.



\bibitem{Reutenauer} C. Reutenauer, Free Lie algebras. London Mathematical
Society Monographs. New Series, 7. Oxford Science Publications. The
Clarendon Press, Oxford University Press, New York, (1993).


\bibitem{Revoy} Philippe Revoy, Alg$\grave{\mbox e}$bres enveloppantes des formes alter$\acute{\mbox e}$es
et des alg$\grave{\mbox e}$bres de Lie, \emph{Journal of Algebra},
\textbf{49}(1977), 342-356.

\bibitem{Shir53} A.I. Shirshov, On the representation of Lie rings in associative rings, \emph{Uspekhi Mat.
Nauk N.S.}, {\bf 8}(1953), No.5(57), 173-175.

\bibitem{Shir1} A.I. Shirshov, On free Lie rings, \emph{ Mat. Sb}., {\bf 45} (1958), 113--122 (in Russian).

\bibitem{Sh62a}A.I. Shirshov, Some algorithmic problem for $\varepsilon$-algebras,
 {\it Sibirsk. Mat. Z.}, \textbf{3}(1962), 132-137.



\bibitem{Sh62b}A.I. Shirshov, Some algorithmic problem for Lie algebras,
{\it Sibirsk. Mat. Z.}, {\bf 3}(2)(1962), 292-296 (in Russian).
English translation: \emph{SIGSAM Bull.},  {\bf 33}(2)(1999), 3-6.

\bibitem{Sh62c}A.I. Shirshov, On the bases of a free Lie algebra, {\it Algebra Logika},
{\bf 1}(1)(1962), 14-19 (in Russian).



\bibitem{Shir3} Selected works of A.I. Shirshov, Eds L.A. Bokut, V. Latyshev, I. Shestakov,
E. Zelmanov, Trs M. Bremner, M. Kochetov, Birkh\"auser, Basel,
Boston, Berlin,  2009.

\bibitem{U95}V.A. Ufnarovski, Combinatorial and asymptotic methods in algebra, Algebra
 VI(57)(1995), 1-196. Springer, Berlin Heidelberg New York.

\bibitem{VG} G. Viennot, Algebras de Lie libres et monoid libres. Bases des
Lie algebres et facrorizations des monoides libres. Lecture Notes in
Mathematics, 691. Berlin, Heidelberg, New York, Springer-Verlag, 124
p. (1978).

\bibitem{Witt} E. Witt, Treue Darstellung Liescher Ringe, \emph{J. reine u. angew. Math.}, \textbf{177}(1937), 152-160.


\end{thebibliography}
\end{document}